\documentclass[11pt]{article}
\usepackage{amssymb,amsmath}
\usepackage{arrows}
\numberwithin{equation}{section}
\allowdisplaybreaks[1]
\newtheorem{theorem}{Theorem}[section]
\newtheorem{corollary}[theorem]{Corollary}
\newtheorem{lemma}[theorem]{Lemma}
\newtheorem{proposition}[theorem]{Proposition}
\newtheorem{conjecture}[theorem]{Conjecture}
\newtheorem{remark}[theorem]{Remark}

\newcommand{\R}{\mathbb{R}}
\newcommand{\C}{\mathbb{C}}

\newcommand{\N}{\mathbb{N}}
\renewcommand{\div}{\mathop{\mathrm{div}}}
\renewcommand{\:}{\thinspace :}

\renewcommand{\d}{\,\mathrm{d}}
\newcommand{\dd}{\mathrm{d}}
\newcommand{\tv}{\mathrm{tv}}
\newcommand{\TS}{\textstyle}
\newcommand{\DS}{\displaystyle}
\newcommand{\app}{\mathrm{app}}
\newcommand{\ran}{\mathrm{ran}}
\newcommand{\weakto}{\rightharpoonup}
\renewcommand{\Re}{\mathop{\mathrm{Re}}}
\renewcommand{\Im}{\mathop{\mathrm{Im}}}
\newcommand{\cA}{\mathcal{A}}
\newcommand{\cE}{\mathcal{E}}
\newcommand{\cL}{\mathcal{L}}
\newcommand{\cM}{\mathcal{M}}
\newcommand{\cO}{\mathcal{O}}
\newcommand{\cS}{\mathcal{S}}
\newcommand{\cT}{\mathcal{T}}
\newcommand{\QED}{\mbox{}\hfill$\Box$}

\begin{document}

\title{Stability and Interaction of Vortices \\
      in Two-Dimensional Viscous Flows}

\author{{\bf Thierry Gallay} \\[1mm] 
Universit\'e de Grenoble I\\
Institut Fourier, UMR CNRS 5582\\
B.P. 74\\
F-38402 Saint-Martin-d'H\`eres, France\\
{\tt Thierry.Gallay@ujf-grenoble.fr}}

\maketitle

\begin{abstract}
The aim of these notes is to present in a comprehensive and relatively
self-contained way some recent developments in the mathematical
analysis of two-dimensional viscous flows. We consider the
incompressible Navier-Stokes equations in the whole plane $\R^2$, and
assume that the initial vorticity is a finite measure. This general
setting includes vortex patches, vortex sheets, and point vortices. We
first prove the existence of a unique global solution, for any value of 
the viscosity parameter, and we investigate its long-time behavior. We
next consider the particular case where the initial flow is a finite
collection of point vortices. In that situation, we show that the
solution behaves, in the vanishing viscosity limit, as a superposition
of Oseen vortices whose centers evolve according to the
Helmholtz-Kirchhoff point vortex system. The proof requires a careful
stability analysis of the Oseen vortices in the large Reynolds number
regime, as well as a precise computation of the deformations of the
vortex cores due to mutual interactions.
\end{abstract}

\tableofcontents

\section{Introduction}
\label{s1}

Although real flows are always three-dimensional, it sometimes happens
that the motion of a fluid is essentially planar in the sense that the
fluid velocity in some distinguished spatial direction is negligible
compared to the velocity in the orthogonal plane.  This situation
often occurs for fluids in thin layers, or for rapidly rotating fluids
where the Coriolis force strongly penalizes displacements along the
axis of rotation. Typical examples are geophysical flows, for which
the geometry of the domain (the atmosphere or the ocean) and the
effect of the Earth's rotation concur to make a two-dimensional
approximation accurate and efficient \cite{CDGG}.
 
From a mathematical point of view, planar flows are substantially
easier to study than three-dimensional ones. For instance, it is known
since the pioneering work of J. Leray \cite{Le1,Le2} that the
two-dimensional incompressible Navier-Stokes equations are globally
well-posed in the energy space, whereas global well-posedness is still
an open problem in the three-dimensional case, no matter which
function space is used \cite{Tao}. The situation is essentially the
same for the incompressible Euler equations, which describe the motion
of inviscid fluids \cite{MB,Con}.  However, having global solutions at
hand does not mean that we fully understand the dynamics of the
system. As a matter of fact, we are not able to establish on a
rigorous basis the phenomenological laws of two-dimensional freely
decaying turbulence, and the stability properties of 2D boundary
layers in the high Reynolds number regime are not fully understood.

In these notes, we consider the idealized situation of an
incompressible viscous fluid filling the whole plane $\R^2$ and
evolving freely without exterior forcing. Following the approach of
Helmholtz \cite{He}, we use the vorticity formulation of the problem,
which is more appropriate to investigate the qualitative behavior of
the solutions. Our goal is to understand the stability properties and
the interactions of localized vortical structures at high Reynolds
numbers. This question is important because carefully controlled
experiments \cite{Cou,RWG}, as well as numerical simulations of
two-dimensional freely decaying turbulence \cite{McW,BMMPW}, suggest
that vortex interactions, and especially vortex mergers, play a
crucial role in the long-time dynamics of viscous planar flows, and
are responsible in particular for the inverse energy cascade. Although
nonperturbative phenomena such as vortex mergers may be very hard to
describe mathematically \cite{MDL,DV}, we shall see that vortex
interactions can be rigorously studied in the asymptotic regime where
the distances between the vortex centers are much larger than the
typical size of the vortex cores.

We begin our analysis of the two-dimensional viscous vorticity
equation in Section~\ref{s2}. We first recall standard estimates for
the two-dimensional Biot-Savart law, which expresses the velocity
field in terms of the vorticity distribution, and we enumerate in
Section~\ref{ss2.1} some general properties of the vorticity
equation, including conservation laws, Lyapunov functions, and scaling
invariance. In Section~\ref{ss2.2}, we introduce the space of
finite measures, which allows to consider nonsmooth flows such as
vortex patches, vortex sheets, or point vortices. It is a remarkable
fact that the vorticity equation is globally well-posed in such a
large function space, for any value of the viscosity parameter, see
\cite{GMO,GG} and Theorem~\ref{thm2} below. Although a complete proof
of that result is beyond the scope of these notes, we show in
Section~\ref{ss2.3} that the classical approach of Fujita and Kato
\cite{FK} applies to our problem and yields global existence and
uniqueness of the solution provided that the atomic part of the
initial vorticity distribution is sufficiently small. For larger
initial data, existence of a global solution can be proved by an
approximation scheme \cite{Co,GMO,Ka2} but additional arguments are
needed to establish uniqueness \cite{GG}.

Section~\ref{s3} collects a few results which describe the behavior of
global solutions of the vorticity equation in $L^1(\R^2)$.  In
particular we prove convergence as $t \to \infty$ to a family of
self-similar solutions called {\em Lamb-Oseen vortices}, see
\cite{GW2} and Theorem~\ref{thm3} below.  To establish these results,
we use accurate estimates on the fundamental solution of
convection-diffusion equations, which were obtained by Osada \cite{Os}
and are reproduced in Section~\ref{ss3.1}. Another fundamental tool is
a transformation into self-similar variables, which compactifies the
trajectories of the system and allows to consider $\omega$-limit sets,
see Section~\ref{ss3.2}. Using a pair of Lyapunov functions, one of
which is only defined for positive solutions, we establish a
``Liouville Theorem'', which characterizes all complete trajectories
of the vorticity equation in $L^1(\R^2)$ that are relatively compact
in the self-similar variables. The conclusion is that all these
trajectories necessarily coincide with Oseen vortices, see
Proposition~\ref{liouville}. In Section~\ref{ss3.3}, we show that
Liouville's theorem implies Theorem~\ref{thm3}, and proves at the same
time that the vorticity equation has a unique solution when the
initial flow is a point vortex of arbitrary circulation. This is
an important particular case of Theorem~\ref{thm2}, which cannot be
established by a standard application of Gronwall's lemma.

In Section~\ref{s4} we investigate in some detail the stability
properties of the Oseen vortices, which are steady states of the
vorticity equation in the self-similar variables.  We introduce in
Section~\ref{ss4.1} an appropriate weighted space for the admissible
perturbations, and we show that the linearized operator at Oseen's
vortex has a remarkable structure in that space\: it is the sum of a
self-adjoint operator, which is essentially the harmonic oscillator,
and a skew-symmetric relatively compact perturbation, which is
multiplied by the circulation of the vortex. This structure almost
immediately implies that Oseen vortices are stable equilibria of the
rescaled vorticity equation, and that the size of the local basin of
attraction is uniform in the circulation parameter. This is in sharp
contrast with many classical examples in fluid mechanics, such as the
Poiseuille flow or the Couette-Taylor flow, for which hydrodynamic
instabilities are known to occur when the Reynolds number becomes
large \cite{DR,TE}. In the case of Oseen vortices, we show in
Section~\ref{ss4.2} that a rapid rotation (i.e., a large circulation
number) has a {\em stabilizing effect} on the vortex\: the size of the
local basin of attraction increases, and the non-radially symmetric
perturbations have a faster decay as $t \to \infty$. These empirically
known facts can be rigorously established, although optimal spectral 
and pseudospectral estimates are not available yet.

In the final section of these notes, we consider the particular
situation where the initial vorticity is a superposition of $N$ Dirac
masses (or point vortices). The corresponding solution of the
two-dimensional vorticity equation is called the {\em viscous
  $N$-vortex solution}, and the goal of Section~\ref{s5} is to
investigate its behavior in the vanishing viscosity limit.  Our main
result, which is stated in Section~\ref{ss5.1}, asserts that the
viscous $N$-vortex solution is nicely approximated by a superposition
of Oseen vortices whose centers evolve according to the
Helmholtz-Kirchhoff point vortex dynamics \cite{MP2,Ne}. This
approximation is accurate as long as the distance between the vortex
centers remains much larger than the typical size of the vortex cores,
which increases through diffusion. In Section~\ref{ss5.2}, we
decompose the viscous $N$-vortex solution into a sum of Gaussian
vortex patches, and we introduce appropriate self-similar variables
which allow us to formulate a stronger version of our result, 
taking into account the deformation of the vortices due to mutual
interactions. The proof involves many technical issues which cannot be
addressed here, but we sketch the main arguments in
Section~\ref{ss5.3} and refer the interested reader to \cite{Ga} for
more details. In particular, we show in Section~\ref{ss5.3} how to
systematically construct an asymptotic expansion of the viscous
$N$-vortex solution, and we briefly indicate how the error terms can
be controlled once a sufficiently accurate approximation is obtained.

The content of the present notes is strongly biased toward the
scientific interests of the author, and does not provide a
comprehensive survey of all important questions in two-dimensional
fluid dynamics. We chose to focus on self-similar vortices, but other
types of flows such as vortex patches or vortex sheets also lead to
interesting and challenging problems, especially in the vanishing
viscosity limit. Also, we should keep in mind that all real fluids are
contained in domains with boundaries, so a comprehensive discussion of
two-dimensional fluid mechanics should certainly include a description
of the flow near the boundary, a question that is totally eluded
here. On the other hand, we do not claim for much originality in these
notes\: all results collected here have already been published
elsewhere, although they were never presented together in a unified
way. In particular, most of the results of Sections~\ref{s3} and
\ref{s4} were obtained in collaboration with C.E. Wayne
\cite{GW1,GW2}, and the content of Section~\ref{s5} is entirely
borrowed from \cite{Ga}. The material presented in Section~\ref{s2} is
rather standard, although our proof of Theorem~\ref{thm1} is perhaps
not explicitly contained in the existing literature. The uniqueness
part in Theorem~\ref{thm2}, which is briefly discussed at the end of
Section~\ref{ss3.3}, was obtained in collaboration with I.~Gallagher
\cite{GG}. Finally, in Section~\ref{ss4.2}, our approach to study the
properties of the linearized operator at Oseen's vortex in the large
circulation regime was developed in a collaboration with I.~Gallagher
and F.~Nier \cite{GGN}.

\section{The Cauchy Problem for the 2D Vorticity Equation}
\label{s2}

We consider the two-dimensional incompressible Navier-Stokes 
equations\:
\begin{equation}\label{NS}
  \left\{\begin{aligned}  
  &\partial_t u(x,t) + (u(x,t)\cdot\nabla)u(x,t) 
  \,=\, \nu \Delta u(x,t) - \frac{1}{\rho}\nabla p(x,t)~, \\
  &\div u(x,t) \,=\, 0~, 
  \end{aligned}\right.
\end{equation}
where $x \in \R^2$ is the space variable and $t \ge 0$ is the time
variable. The unknown functions are the velocity field $u(x,t) = 
(u_1(x,t),u_2(x,t)) \in \R^2$, which represents the speed of 
a fluid particle at point $x$ and time $t$, and the pressure field 
$p(x,t) \in \R$. Eq.~\eqref{NS} contains two physical parameters, the 
kinematic viscosity $\nu > 0$ and the fluid density $\rho > 0$, 
which are both assumed to be constant. 

Equivalently, the motion of a planar fluid can be described by 
the {\em vorticity} field\:
\[
  \omega(x,t) \,=\, \partial_1 u_2(x,t) - \partial_2 u_1(x,t)~,
\]
which represents the angular rotation of the fluid particles. 
If we take the two-dimensional curl of the first equation in \eqref{NS}, 
we obtain the evolution equation
\begin{equation}\label{omeq}
  \partial_t \omega(x,t) +  u(x,t)\cdot\nabla \omega(x,t) \,=\, 
  \nu \Delta \omega(x,t)~,
\end{equation}
which is the starting point of our analysis. Note that $u \cdot \nabla
\omega = \div(u\omega)$ because $\div u = 0$. In the case of a
perfect fluid ($\nu = 0$), Eq.~\eqref{omeq} simply means that the
vorticity is advected by the velocity field $u(x,t)$ like a material
particle. For real fluids ($\nu > 0$), the vorticity also diffuses at
a rate given by the kinematic viscosity.
 
The vorticity equation \eqref{omeq} is definitely simpler than the
original Navier-Stokes system \eqref{NS}, but it still contains 
the velocity field $u(x,t)$. To make \eqref{omeq} independent of 
\eqref{NS}, it is possible to express $u$ in terms of $\omega$ by 
solving the elliptic system
\begin{equation}\label{uomega}
  \partial_1 u_1 + \partial_2 u_2 \,=\, 0~, \qquad
  \partial_1 u_2 - \partial_2 u_1 \,=\, \omega~. 
\end{equation}
If $\omega$ decays sufficiently fast at infinity (see Lemma~\ref{lemBS}
below), the solution is given by the two-dimensional {\em Biot-Savart} 
formula\:
\begin{equation}\label{BS}
  u(x,t) \,=\, \frac{1}{2\pi}\int_{\R^2} \frac{(x-y)^\perp}{|x-y|^2}
  \,\omega(y,t) \d y~,
\end{equation}
where $x^\perp = (-x_2,x_1)$ and $|x|^2 = x_1^2 + x_2^2$ if 
$x = (x_1,x_2) \in \R^2$. The vorticity equation \eqref{omeq}, 
supplemented with the Biot-Savart law \eqref{BS}, is now formally
equivalent to the Navier-Stokes equations \eqref{NS}. Once a 
solution to \eqref{omeq}, \eqref{BS} is found, the pressure 
field can be recovered up to an additive constant by solving, 
for each $t > 0$, the Poisson equation
\[
  -\Delta p \,=\, \rho \div ((u \cdot \nabla) u)~,
  \qquad p(x,t) \to 0 \quad\hbox{as}\quad |x| \to \infty~.
\]

The following standard result shows that the Biot-Savart formula
\eqref{BS} is well defined if the vorticity $\omega$ lies in 
$L^p(\R^2)$ for some $p \in (1,2)$. 

\begin{lemma}\label{lemBS}
Assume that $\omega \in L^p(\R^2)$ for some $p \in (1,2)$. Then 
the velocity field $u$ defined (for almost every $x \in \R^2$) 
by the Biot-Savart formula \eqref{BS} satisfies\:\\
i) $u \in L^q(\R^2)$, where $q \in (2,\infty)$ is such that
$\frac{1}{q} = \frac{1}{p} - \frac{1}{2}$. Moreover
\begin{equation}\label{BSi}
  \|u\|_{L^q(\R^2)} \,\le\, C \|\omega\|_{L^p(\R^2)}~.
\end{equation}
ii) $\nabla u \in L^p(\R^2)$, and \eqref{uomega} holds in $L^p(\R^2)$. 
Moreover
\begin{equation}\label{BSii}
  \|\nabla u\|_{L^p(\R^2)} \,\le\, C \|\omega\|_{L^p(\R^2)}~. 
\end{equation}
\end{lemma}

\noindent{\bf Proof.} 
According to \eqref{BS} we have $u = K * \omega$, where 
\[
  K(x) \,=\, \frac{1}{2\pi}\,\frac{x^\perp}{|x|^2}~, \qquad 
  x \in \R^2\setminus\{0\}~, 
\]
is the Biot-Savart kernel. Since $K$ belongs to the weak 
$L^2$ space $L^{2,\infty}(\R^2)$, the Hardy-Littlewood-Sobolev
(or weak Young) inequality \cite[Section 4.3]{LL} shows that 
\[
  \|u\|_{L^q(\R^2)} \,=\, \|K*\omega\|_{L^q(\R^2)} \,\le\, 
  \|K\|_{L^{2,\infty}(\R^2)}\|\omega\|_{L^p(\R^2)}~,
\]
if $1 < p < 2 < \infty$ and $\frac{1}{q} = \frac{1}{p} - \frac{1}{2}$.
This proves \eqref{BSi}. Moreover $\nabla u = \nabla K * \omega$, 
where $\nabla K$ is an integral kernel of Calder\'on-Zygmund type
which is associated to a bounded Fourier multiplier. It follows 
\cite[Section I.5]{St} that $K$ defines a bounded linear operator
in $L^p(\R^2)$ for $p \in (1,2)$, as asserted in \eqref{BSii}.  
\QED

\begin{remark}\label{RemBS}
If $\omega \in L^1(\R^2)$, the velocity field $u = K*\omega$ 
belongs to $L^{2,\infty}(\R^2)$, but $u \notin L^2(\R^2)$ in general, 
see \cite{GMO}. If $\omega \in L^2(\R^2)$, one can solve \eqref{uomega}
directly by a simple calculation in Fourier space, which shows that 
$u$ belongs to the homogeneous Sobolev space $\dot H^1(\R^2)$, but 
$u \notin L^\infty(\R^2)$ in general. Finally, if $\omega \in L^p(\R^2)$ 
for some $p \in (2,\infty)$, the solution of \eqref{uomega} is no 
longer given by the Biot-Savart law \eqref{BS}, but by a modified 
formula of the form
\[
  u(x) \,=\, \frac{1}{2\pi}\int_{\R^2} \left(\frac{(x-y)^\perp}{|x-y|^2}
  + \frac{y^\perp}{|y|^2}\right)\,\omega(y) \d y~.
\]
In that case, $u(x)$ may grow like $|x|^{1-\frac2p}$ as $|x| \to 
\infty$, and the solution of \eqref{uomega} is only defined up 
to an additive constant.
\end{remark}

\subsection{General Properties of the Vorticity Equation}\label{ss2.1}

In this section, we list a few important properties of the vorticity
equation \eqref{omeq}, including conservation laws, Lyapunov functions, 
and scaling invariance. The presentation is formal in the sense 
that we assume here that we are given a solution of \eqref{omeq} with
the appropriate smoothness and integrability properties. However, 
all calculations can be justified once suitable existence theorems 
have been established \cite{MB}. 

\subsubsection{Conservations Laws} Let $\omega(x,t)$ be a solution 
of \eqref{omeq} with initial data $\omega_0(x) = \omega(x,0)$. 
If $\omega_0 \in L^1(\R^2)$, then $x \mapsto \omega(x,t)$ is integrable 
for all $t \ge 0$ and the {\em total circulation} is conserved\:
\begin{equation}\label{totcir}
  \int_{\R^2}\omega(x,t)\d x \,=\, \int_{\R^2}\omega_0(x)\d x~. 
\end{equation}
If moreover $|x|\omega_0 \in L^1(\R^2)$, then the first order moments
of $\omega(x,t)$ are also conserved\:
\begin{equation}\label{firstmom}
  \int_{\R^2}x_i \omega(x,t)\d x \,=\, \int_{\R^2}x_i \omega_0(x)\d x~,
  \qquad i = 1,2~. 
\end{equation}
Finally, if in addition $|x|^2 \omega_0 \in L^1(\R^2)$, then
the symmetric second order moment of $\omega(x,t)$ satisfies\:
\begin{equation}\label{secondmom}
  \int_{\R^2}|x|^2 \omega(x,t)\d x \,=\, \int_{\R^2}|x|^2 \omega_0(x)\d x
  + 4\nu t \int_{\R^2}\omega_0(x)\d x~. 
\end{equation}
Thus the symmetric moment is conserved for perfect fluids ($\nu = 0$),
and for viscous fluids ($\nu > 0$) if the total circulation vanishes.
These elementary properties are easily established by direct
calculations, using \eqref{omeq} and the Biot-Savart law $u = K *
\omega$. In particular, to prove \eqref{firstmom}, one uses the fact
that
\[
   \int_{\R^2}u_i(x) \omega(x) \d x \,=\, \int_{\R^2 \times \R^2}
   K_i(x-y) \omega(y)\omega(x)\d x \d y \,=\, 0~, \qquad
   i = 1,2~,
\]
because the Biot-Savart kernel $K$ is odd. Similarly, the proof of 
\eqref{secondmom} relies on the identity
\begin{equation}\label{secmomid}
   \int_{\R^2}(x\cdot u(x)) \omega(x) \d x \,=\, \frac12 
   \int_{\R^2 \times \R^2} \Bigl((x-y) \cdot K(x-y)\Bigr) \omega(y)
   \omega(x)\d x \d y \,=\, 0~,
\end{equation}
which holds because $K$ is odd and $x \cdot K(x) = 0$ for all 
$x \in \R^2$. 

\begin{remark}\label{momrem}
By Green's theorem, the total circulation satisfies 
\[
  \gamma \,:=\, \int_{\R^2}\omega\d x \,=\, \lim_{R \to \infty} 
  \oint_{|x|=R}(u_1(x)\d x_1 + u_2(x)\d x_2)~,
\]
hence $\gamma$ represents the circulation of the velocity field
at infinity. This quantity is conserved due to the structure 
of Eq.~\eqref{omeq}. In contrast, the conservation of the first-order 
moments and the simple evolution law for the symmetric second order 
moment are related to translation and rotation invariance, and 
would not hold if the vorticity equation was considered in a 
nontrivial domain $\Omega \subset \R^2$. 
\end{remark}

\subsubsection{Lyapunov Functions} If $\omega_0 \in L^p(\R^2)$ 
for some $p \in [1,\infty]$, the solution of \eqref{omeq} with initial 
data $\omega_0$ lies in $L^p(\R^2)$ for all times and 
\begin{equation}\label{Lpdecay}
  \|\omega(\cdot,t)\|_{L^p(\R^2)} \,\le\,  \|\omega_0\|_{L^p(\R^2)}~,
  \quad \hbox{for all }t \ge 0~.
\end{equation}
This again follows from the fact that the advection field $u$ in
\eqref{omeq} is divergence-free.  If the vorticity equation 
is considered in a domain $\Omega \subset  \R^2$, and if 
one assumes as usual that $u = 0$ on $\partial\Omega$
(no-slip boundary condition), then vorticity is created in the 
boundary layer near $\partial\Omega$ and \eqref{Lpdecay}
does not hold. 

If $u_0 = K*\omega_0 \in L^2(\R^2)$, the solution $u(x,t)$ 
of the Navier-Stokes equation \eqref{NS} with initial data 
$u_0$ satisfies $u(\cdot,t) \in L^2(\R^2)$ for all times, 
and the {\em kinetic energy} 
\[
  E(t) \,=\, \frac12 \int_{\R^2} |u(x,t)|^2\d x~, 
 \]
is a nonincreasing function of time\:
\[
  E'(t) \,=\, -\nu\int_{\R^2} |\nabla u(x,t)|^2\d x \,\le\, 0~.
\]
The kinetic energy is the most famous Lyapunov function in fluid
mechanics, and often the only one that is available. In particular, it
plays a crucial role in the construction of global weak solutions of
the three-dimensional Navier-Stokes equations \cite{Le2}. Somewhat
paradoxically, in unbounded two-dimensional domains, this quantity is
not always useful because it is finite only if the total circulation
vanishes. Indeed, we have the following elementary result\:

\begin{lemma}\label{kinerg}
If $u \in L^2(\R^2)^2$ satisfies $\omega = \partial_1 u_2 - 
\partial_2 u_1 \in L^1(\R^2)$, then $\int_{\R^2} \omega \d x = 0$. 
\end{lemma}

\noindent{\bf Proof.}
By assumption the Fourier transform $\hat \omega(k) 
= ik_1 \hat u_2(k) - ik_2 \hat u_1(k)$ is a continuous 
function of $k \in \R^2$, hence
\[
  \int_{\R^2}\omega \d x \,=\, \hat \omega(0) \,=\, 
  \lim_{\epsilon\to 0} \frac{1}{\pi\epsilon^2}\int_{|k|\le\epsilon} 
  \hat\omega(k)\d k~.
\]
But
\[
  \frac{1}{\pi\epsilon^2}\int_{|k|\le\epsilon} |\hat\omega(k)|
  \d k \,\le\, \frac{1}{\pi\epsilon^2}\int_{|k|\le\epsilon}  
  |k| |\hat u(k)|\d k \,\le\, \frac{1}{\sqrt{2\pi}} 
  \Bigl(\int_{|k|\le\epsilon} |\hat u(k)|^2 \d k\Bigr)^{1/2}~,
\]
thus taking the limit $\epsilon \to 0$ we obtain $\int_{\R^2} \omega 
\d x = 0$.  
\QED

As a substitute for the kinetic energy, one can consider 
the {\em pseudo-energy}
\[
  \cE_d(t) \,=\, \frac{1}{4\pi} \int_{\R^2\times\R^2}
  \log\frac{d}{|x-y|}\,\omega(x,t)\omega(y,t)\d x\d y~, 
\]
where $d > 0$ is an arbitrary length scale. This quantity is well
defined under mild integrability assumptions on $\omega$, and can be
finite even if the total circulation is nonzero. Morever, $\cE_d(t)$
is a Lyapunov function\:
\[
  \cE_d'(t) \,=\, -\nu \int_{\R^2}\omega(x,t)^2\d x \,\le\, 0~,
\]
and it is straightforward to verify that $\cE_d(t) = E(t)$ 
if the velocity field associated to $\omega(x,t)$ is square 
integrable. In this particular situation, $\cE_d(t)$ does not depend 
on the parameter $d$.

\subsubsection{Scaling Invariance}

In the whole plane $\R^2$, the Navier-Stokes equations \eqref{NS}
and the vorticity equation \eqref{omeq} are invariant under the 
rescaling
\begin{equation}\label{scaling}
  u(x,t) \,\mapsto \lambda u(\lambda x,\lambda^2 t)~, \qquad
  \omega(x,t) \,\mapsto \lambda^2 \omega(\lambda x,\lambda^2 t)~, 
\end{equation}
for any $\lambda > 0$. To solve the Cauchy problem, it is rather
natural to use space-time norms which are invariant under the 
transformation \eqref{scaling}. The corresponding function 
spaces are called {\em critical spaces} (with respect 
to the scaling). A typical example is the {\em energy space}
\begin{equation}\label{enspace}
  u \in X = C^0_b(\R_+,L^2(\R^2))~, \qquad \|u\|_X \,=\, 
  \sup_{t \ge 0} \|u(\cdot,t)\|_{L^2}~, 
\end{equation}
where $C^0_b(\R_+,L^2(\R^2))$ is the space of all bounded and
continuous maps from $\R_+ = [0,\infty)$ into $L^2(\R^2)$, and
$u(\cdot,t)$ denotes the map $x \mapsto u(x,t)$. Alternatively, to
include solutions with nonzero total circulation, one can assume that
the vorticity is integrable and use the space
\begin{equation}\label{omspace}
  \omega \in Y = C^0_b(\R_+,L^1(\R^2))~, \qquad \|\omega\|_Y \,=\, 
  \sup_{t \ge 0}\|\omega(\cdot,t)\|_{L^1}~.  
\end{equation}

A well-known ``metatheorem'' in nonlinear PDE's asserts that, for a 
nonlinear partial differential equation with scaling invariance, 
the critical spaces are the largest ones (in terms of local regularity 
of the solutions) in which we can hope that the Cauchy problem is 
locally well-posed. In fact, there may be other obstructions to 
local well-posedness, but according to the above claim it is certainly 
not reasonable to try to solve the Cauchy problem in function spaces 
that are ``strictly larger'' than the critical ones. 

As far as the Navier-Stokes equation is concerned, we know since the
work of Leray \cite{Le1} that Eq.~\eqref{NS} is globally well-posed in
the energy space \eqref{enspace}, and a similar result holds for the
vorticity equation \eqref{omeq} in the critical space \eqref{omspace},
see \cite{BA,Br}. In the rest of this section, we discuss the Cauchy
problem for \eqref{omeq} in the space of finite measures, which 
is a natural extension of $L^1(\R^2)$. As we shall see in 
Section~\ref{s5}, this generalization is essential to study the 
interaction of vortices in the vanishing viscosity limit. 

\subsection{The Space of Finite Measures}\label{ss2.2}

Let $\cM(\R^2)$ denote the space of all real-valued 
Radon measures on $\R^2$, equipped with the total 
variation norm\:
\begin{equation}\label{tvdef}
  \|\mu\|_\tv \,=\, \sup\Bigl\{\int_{\R^2}\phi\d\mu ~\Big|~
  \phi \in C_0(\R^2)\,,~\|\phi\|_{L^\infty} \le 1\Bigr\}~.
\end{equation}
Here $C_0(\R^2)$ is the space of all continuous functions $\phi : \R^2
\to \R$ such that $\phi(x) \to 0$ as $|x| \to \infty$.  It is well
known that $\cM(\R^2)$ equipped with the total variation norm is a
Banach space, which contains $L^1(\R^2)$ as a closed subspace (if we
identify any absolutely continuous measure with its density
function). In particular, if $\omega \in L^1(\R^2)$, then
$\|\omega\|_\tv = \|\omega\|_{L^1}$. The total variation norm is scale
invariant in the sense that $\|\mu_\lambda\|_\tv = \|\mu\|_\tv$ for
all $\lambda > 0$, where $\mu_\lambda$ is the rescaled measure defined
by
\[
  \int_{\R^2}\phi\d\mu_\lambda \,=\, \int_{\R^2}\phi(\lambda^{-1}\cdot)
  \d\mu~, \qquad \hbox{for all } \phi \in C_0(\R^2)~.
\]
As is clear from the definition \eqref{tvdef}, the space $\cM(\R^2)$ 
is the tolopogical dual of $C_0(\R^2)$. By the Banach-Alaoglu 
theorem, it follows that the unit ball in $\cM(\R^2)$ is a (sequentially) 
compact set for the weak convergence defined by\:
\[
  \mu_n \,\xrightharpoonup[n \to \infty]{\hbox to 6mm{}}\, \mu
  \qquad \hbox{if}\qquad \int_{\R^2}\phi\d\mu_n 
  \,\xrightarrow[n \to \infty]{}\,
  \int_{\R^2}\phi\d\mu \quad \hbox{for all }\phi \in C_0(\R^2)~. 
\]  

Any $\mu \in \cM(\R^2)$ can be decomposed in a unique way as 
$\mu = \mu_{ac} + \mu_s$, where $\mu_{ac} \in \cM(\R^2)$ is 
absolutely continuous with respect to Lebesgue's measure on 
$\R^2$, and $\mu_s$ is singular with respect to Lebesgue's 
measure \cite{Ru}. By the Radon-Nikodym theorem, there exists
a unique $\omega \in L^1(\R^2)$ such that $\dd\mu_{ac} = \omega
\d x$. Furthermore, the singular part $\mu_s$ can be decomposed
as $\mu_s = \mu _{sc} + \mu_{pp}$, where $\mu_{pp}$ is the restriction 
of $\mu$ to the atomic set $\Sigma = \{x \in \R^2 \,|\, \mu(\{x\}) 
\neq 0\}$. Since $\mu$ is a finite measure, this set is at most
countable, so that 
\begin{equation}\label{muppdecomp}
  \mu_{pp} \,=\, \sum_{i=1}^\infty \gamma_i \delta_{x_i}~,
\end{equation}
for some $\gamma_i \in \R$ and some $x_i \in \R^2$. By construction, 
the singularly continuous part $\mu_{sc}$ has no atoms, yet is 
concentrated on a set of zero Lebesgue measure. 

To summarize, we have $\mu = \mu_{ac} + \mu_{sc} + \mu_{pp}$, with 
$\mu_{ac} \perp \mu_{sc} \perp \mu_{pp}$. This notation means that 
the absolutely continuous, singularly continuous, and pure point 
parts of $\mu$ are concentrated on pairwise disjoint sets. 
It follows in particular that
\[
  \|\mu\|_\tv \,=\, \|\mu_{ac}\|_\tv + \|\mu_{sc}\|_\tv + 
  \|\mu_{pp}\|_\tv \,=\, \|\omega\|_{L^1} + \|\mu_{sc}\|_\tv + 
  {\TS\sum\limits_{i=1}^\infty |\gamma_i|}~.  
\]

If we now suppose that the measure $\mu \in \cM(\R^2)$ is the 
vorticity of a two-dimensional flow, there is a nice correspondence
between the abstract decomposition result presented above and the 
following standard catalogue of nonsmooth flows\:

\begin{itemize}

\item {\em Vortex Patches.} A vortex patch can be defined as 
a flow for which the vorticity $\mu$ is a piecewise smooth function. 
For instance, $\mu$ can be the characteristic function 
of a smooth bounded domain $\Omega \subset \R^2$. In that case 
$\mu = \mu_{ac} \in L^1(\R^2) \cap L^\infty(\R^2)$, and the velocity 
field $u = K * \mu$ is Lipschitz continuous if the boundary $\partial
\Omega$ is sufficiently smooth. 

\item {\em Vortex Sheets.} In contrast, vortex sheets are characterized
by discontinuities of the velocity field. As a typical example, 
we mention the case where the vorticity $\mu$ is the Euclidean 
measure supported by a smooth closed curve $\gamma \subset \R^2$.
In that case $\mu = \mu_{sc}$, and the tangential component 
of the velocity field $u = K*\mu$ is discontinuous on the curve 
$\gamma$. 

\item {\em Point Vortices.} This is the case where the vorticity 
$\mu$ is a collection of Dirac masses. Here $\mu = \mu_{pp}$, 
and the velocity field $u = K*\mu$ is unbounded near the center
of each vortex.

\end{itemize}

\subsection{The Cauchy Problem in $\cM(\R^2)$}\label{ss2.3}

The aim of this section is to study the Cauchy problem for 
the two-dimensional vorticity equation \eqref{omeq} in the 
space of finite measures. We start with a preliminary result that 
is relatively easy to prove. 

\begin{theorem}\label{thm1} {\rm \cite{Co,GMO,Ka2}}\\
There exists a universal constant $C_0 > 0$ such that, if the initial 
vorticity $\mu\in \cM(\R^2)$ satisfies $\|\mu_{pp}\|_\tv \le C_0\,\nu$, 
then the vorticity equation \eqref{omeq} has a unique global solution
\[
  \omega \in C^0((0,\infty),L^1(\R^2) \cap L^\infty(\R^2))
\]
such that $\|\omega(\cdot,t)\|_{L^1} \le \|\mu\|_\tv$ for all 
$t > 0$, and $\omega(\cdot,t) \weakto \mu$ as $t \to 0$.
\end{theorem}

Existence of global solutions to \eqref{omeq} with initial data in
$\cM(\R^2)$ was first proved by G.-H. Cottet \cite{Co}, by Y. Giga,
T. Miyakawa and H. Osada \cite{GMO}, and by T. Kato \cite{Ka2}.  In
these works, the strategy is\: i) to approximate the initial measure
by a smooth vorticity distribution; ii) to construct a global
solution, with good a priori estimates, by applying classical
existence results for smooth initial data; iii) to obtain a solution
of the original problem by extracting an appropriate subsequence as
the regularization parameter converges to zero. This yields the
existence part of Theorem~\ref{thm1} without any smallness condition
on the initial data. The restriction $\|\mu_{pp}\|_\tv \le C_0\,\nu$
arises when one tries to prove uniqueness of the solution by a
standard application of Gronwall's lemma. As we shall see later, this
smallness condition is purely technical and can be completely relaxed,
see Theorem~\ref{thm2}.

In the rest of this section, we present a more direct proof of
Theorem~\ref{thm1}, which is nevertheless inspired by \cite{GMO}.
Without loss of generality, we assume from now on that the kinematic
viscosity is equal to $1$. Given initial data $\mu \in \cM(\R^2)$,
we consider the integral equation associated with \eqref{omeq}\:
\begin{equation}\label{omin}
  \omega(t) \,=\, e^{t\Delta}\mu - \int_0^t \div \Bigl(
  e^{(t-s)\Delta} \,u(s)\omega(s)\Bigr)\d s~, \qquad t > 0~,
\end{equation}
where $\omega(t) = \omega(\cdot,t)$, $u(t) = u(\cdot,t)$, and 
$e^{t\Delta}$ denotes the heat semigroup defined by
\begin{equation}\label{heat}
  (e^{t\Delta}\mu)(x) \,=\, \frac{1}{4\pi t}\int_{\R^2}
  e^{-|x-y|^2/(4t)}\d \mu_y~, \qquad t > 0~, \quad x \in \R^2~.
\end{equation}
Our goal is to solve the integral equation \eqref{omin} by a fixed 
point argument in an appropriate function space. This is a classical
method which, in the context of the Navier-Stokes equations, goes 
back to the work of Kato and Fujita \cite{FK}. Solutions of the integral 
equation \eqref{omin} are often called {\em mild} solutions of 
the vorticity equation \eqref{omeq}. It is rather straightforward 
to show that, if $\omega \in C^0((0,T),L^p(\R^2))$ is a mild 
solution of \eqref{omeq} with $p \ge 4/3$, then $\omega(x,t)$ is 
smooth and satisfies \eqref{omeq} in the classical sense. Thus there
is no loss of generality in considering \eqref{omin} instead of 
\eqref{omeq}. 

\subsubsection{Heat Kernel Estimates}

The following lemma plays a crucial role in the proof of 
Theorem~\ref{thm1}\:

\begin{lemma}\label{heatlem} Let $\mu \in \cM(\R^2)$.\\
a) For $1 \le p \le \infty$ and $t > 0$, we have
\begin{equation}\label{heat1}
  \|e^{t\Delta}\mu\|_{L^p} \,\le\, \frac{1}{(4\pi t)^{1-\frac1p}}
  \,\|\mu\|_\tv~, \qquad 
  \|\nabla e^{t\Delta}\mu\|_{L^p} \,\le\, \frac{C}{t^{\frac32 -\frac1p}}
  \,\|\mu\|_\tv~.
\end{equation}
b) For $1 < p \le \infty$, we have
\begin{equation}\label{heat2}
  L_p(\mu) \,:=\, \limsup_{t \to 0} \,(4\pi t)^{1-\frac1p} 
  \|e^{t\Delta}\mu\|_{L^p} \,\le\, \|\mu_{pp}\|_\tv~.
\end{equation}
\end{lemma}

\noindent{\bf Proof.}
Estimates \eqref{heat1} are easily obtained from \eqref{heat} if 
$p = 1$ or $p = \infty$, and the general case follows by interpolation.
Thus we concentrate on estimate \eqref{heat2}, which we establish 
using a simplified version of the argument given by Giga, Miyakawa, 
and Osada \cite{GMO}. First, if we decompose $\mu = \mu_{ac} + \mu_{sc} 
+ \mu_{pp}$, we have
\[
  L_p(\mu) \,\le\, L_p(\mu_{ac}) + L_p(\mu_{sc}) + L_p(\mu_{pp}) 
  \,\le\, L_p(\mu_{ac}) + L_p(\mu_{sc}) + \|\mu_{pp}\|_\tv~. 
\]
Thus, to prove \eqref{heat2}, it is sufficient to show that 
$L_p(\mu) = 0$ for all $\mu \in \cM(\R^2)$ such that $\mu_{pp} = 0$. 
Furthermore, since
\[
  L_p(\mu) \,\le\, L_1(\mu)^{1/p} L_\infty(\mu)^{1-1/p} \,\le\,  
  \|\mu\|_\tv^{1/p} L_\infty(\mu)^{1-1/p}~,
\]
it is sufficient to consider the case where $p = \infty$. 

Assume now that $\mu \in \cM(\R^2)$ satisfies $\mu_{pp} = 0$, 
and fix $\epsilon > 0$. We claim that there exists $\delta > 0$ such 
that
\begin{equation}\label{Bdel}
  \sup_{x \in \R^2} |\mu|(B(x,\delta)) \,\le\, \epsilon~, \quad
  \hbox{where} \quad B(x,\delta) \,=\, \{y \in \R^2\,|\, 
  |y-x| \le \delta\}~.
\end{equation}
Here and in what follows, we denote by $|\mu| \in \cM(\R^2)$ the 
total variation of the measure $\mu \in \cM(\R^2)$. If $\mu^+, 
\mu^-$ are the positive and negative variations of $\mu$, 
then $\mu = \mu^+ - \mu^-$ (Jordan's decomposition) and 
$|\mu| = \mu^+ + \mu^-$ \cite{Ru}. To prove \eqref{Bdel}, assume
on the contrary that there exist $\delta_n > 0$ and $x_n \in \R^2$ 
such that $\delta_n \to 0$ as $n \to \infty$ and $|\mu|(B(x_n,\delta_n)) 
> \epsilon$ for all $n \in \N$. Since $|\mu|$ is a finite measure,
the sequence $(x_n)$ is necessarily bounded, thus after extracting
a subsequence we can assume that $x_n \to \tilde x$ as $n \to \infty$, 
for some $\tilde x \in \R^2$. As a consequence, for any $\delta > 0$, 
we have $B(\tilde x,\delta) \supset B(x_n,\delta_n)$ for all sufficiently
large $n$, so that $|\mu|(B(\tilde x,\delta)) \ge |\mu|(B(x_n,\delta_n)) 
> \epsilon$. Taking the limit $\delta \to 0$, we obtain $|\mu(\{\tilde x
\})| \ge \epsilon$, which contradicts the assumption that $\mu_{pp} = 0$. 
This proves \eqref{Bdel}. 

Now, for any $t > 0$, we can take $\bar x(t) \in \R^2$ such that 
$|(e^{t\Delta}\mu)(\bar x(t))|  = \|e^{t\Delta}\mu\|_{L^\infty}$. 
Thus, using \eqref{heat} and \eqref{Bdel}, we obtain
\[
  4\pi t\,\|e^{t\Delta}\mu\|_{L^\infty} \,\le\, \int_{B(\bar x(t),\delta)}
  e^{-\frac{|\bar x(t)-y|^2}{4t}}\d|\mu|_y + \int_{B(\bar x(t),\delta)^c} 
  e^{-\frac{|\bar x(t)-y|^2}{4t}}\d|\mu|_y~.
\]
The first term on the right-hand side is bounded by $\epsilon$ for 
all $t > 0$, and the second one vanishes as $t \to 0$. Thus 
\[
  L_\infty(\mu) \,=\, \limsup_{t \to 0}\,(4\pi t)\,\|e^{t\Delta}
  \mu\|_{L^\infty} \,\le\, \epsilon~.
\]
Since $\epsilon > 0$ was arbitrary, we conclude that $L_\infty(\mu) = 0$
whenever $\mu_{pp} = 0$, which is the desired result. 
\QED

\begin{remark}\label{improvedrem}
Estimate \eqref{heat2} was strengthened by Kato \cite{Ka2} in 
the following way\:
\begin{equation}\label{heat3}
  \lim_{t \to 0} \,(4\pi t)^{1-\frac1p} \|e^{t\Delta}\mu\|_{L^p} 
  \,=\, p^{-1/p} \|\{\gamma_i\}_{i=1}^\infty\|_{\ell^p} \,\le\, 
  {\TS\sum\limits_{i=1}^\infty |\gamma_i|}~, 
\end{equation}
where $\gamma_i$, $x_i$ are as in \eqref{muppdecomp}. In both 
\eqref{heat2}, \eqref{heat3}, the assumption that $p > 1$ is 
crucial\: if $\mu \in \cM(\R^2)$ is a positive measure, then 
$\|e^{t\Delta}\mu\|_{L^1} = \|\mu\|_\tv$ for all $t > 0$, hence 
$L_1(\mu) = \|\mu\|_\tv$. 
\end{remark}

\subsubsection{Fixed Point Argument}

To prove the existence of solutions to the integral equation 
\eqref{omin}, we fix $T > 0$ and introduce the function space
\begin{equation}\label{XTdef}
  X_T \,=\, \Bigl\{\omega \in C^0((0,T],L^{4/3}(\R^2)) \,\Big|\,
  \|\omega\|_{X_T} < \infty \Bigr\}~, 
\end{equation}
equipped with the norm
\[
  \|\omega\|_{X_T} \,=\, \sup_{0 < t \le T} t^{1/4} 
  \|\omega(t)\|_{L^{4/3}}~.
\]
Given $\mu \in \cM(\R^2)$, we denote $\omega_0(t) = e^{t\Delta}\mu$ 
for all $t > 0$. In view of Lemma~\ref{heatlem}, we have $\omega_0 \in 
X_T$ and there exist positive constants $C_1, C_2$ (independent of 
$T$) such that

\smallskip\noindent
~a) $\|\omega_0\|_{X_T} \le C_1 \|\mu\|_\tv$ (for any $T > 0$);

\smallskip\noindent
~b) $\|\omega_0\|_{X_T} \le C_2 \|\mu_{pp}\|_\tv + \epsilon$ (if 
$T > 0$ is sufficiently small, depending on $\mu$,$\epsilon$).

\smallskip\noindent
Here $\epsilon > 0$ is an arbitrarily small positive constant. 

On the other hand, given $\omega \in X_T$, we define $F\omega 
: (0,T] \to L^{4/3}(\R^2)$ by
\begin{equation}\label{Fdef}
  (F\omega)(t) \,=\, \int_0^t \div \Bigl(e^{(t-s)\Delta} \,u(s)
  \omega(s)\Bigr)\d s~, \quad 0 < t \le T~,
\end{equation}
where $u(s) = K*\omega(s)$ is the velocity field obtained via 
the Biot-Savart law \eqref{BS}. Using the second estimate in 
\eqref{heat1}, H\"older's inequality, and estimate \eqref{BSi}, 
we find
\begin{align}\nonumber
  t^{1/4}\|(F\omega)(t)\|_{L^{4/3}} \,&\le\, 
  t^{1/4} \int_0^t \frac{C}{(t-s)^{\frac34}}
  \,\|u(s) \omega(s)\|_{L^1}\d s \\ \nonumber
  \,&\le\, t^{1/4} \int_0^t \frac{C}{(t-s)^{\frac34}}
  \,\|u(s)\|_{L^4} \|\omega(s)\|_{L^{4/3}}\d s \\ \label{Fest}
  \,&\le\, t^{1/4} \int_0^t \frac{C}{(t-s)^{\frac34}}
  \,\|\omega(s)\|_{L^{4/3}}^2\d s \\ \nonumber
  \,&\le\, C \|\omega\|_{X_T}^2 \,t^{1/4} \!\int_0^t \frac{C }{
  (t-s)^{\frac34}s^{\frac12}}\d s \,\le\, C \|\omega\|_{X_T}^2~,
\end{align}
for all $t \in (0,T]$. This shows that $\|F\omega\|_{X_T} \le C_3 
\|\omega\|_{X_T}^2$ for some $C_3 > 0$ (independent of $T$), and it is 
not difficult to verify that $(F\omega)(t)$ depends continuously on 
$t$ in the topology of $L^{4/3}(\R^2)$, so that $F\omega \in X_T$. 
Similarly, one can show that
\begin{equation}\label{LipF}
  \|F\omega - F\tilde\omega\|_{X_T} \,\le\, C_3 (\|\omega\|_{X_T} + 
  \|\tilde\omega\|_{X_T})\|\omega - \tilde\omega\|_{X_T}~,
\end{equation}
for all $\omega, \tilde \omega \in X_T$. 

Now, fix $R > 0$ such that $2C_3 R < 1$, and consider the closed ball
\[
  B \,=\, \{\omega \in X_T\,|\, \|\omega\|_{X_T} \le R\} 
  \,\subset\, X_T~.
\]
As a consequence of the estimates above, if $\|\omega_0\|_{X_T} \le 
R/2$, the map $\omega \mapsto \omega_0 - F\omega$ is a strict 
contraction in $B$, hence has a unique fixed point in $B$. By 
construction, this fixed point $\omega$ is a solution of the 
integral equation \eqref{omin} on $(0,T]$. The condition 
$\|\omega_0\|_{X_T} \le R/2$ is fulfilled if either

\smallskip\noindent
~1) $2C_1 \|\mu\|_\tv \le R$ and $T > 0$ is arbitrary, or

\smallskip\noindent
~2) $2C_2 \|\mu_{pp}\|_\tv < R$ and $T > 0$ is small enough, 
depending on $\mu$. 

\smallskip\noindent
In other words, the proof above shows that equation \eqref{omin} is
{\em globally well-posed} for small initial data $\mu \in \cM(\R^2)$,
and {\em locally well-posed} for large initial data with small atomic
part $\mu_{pp}$. However, if $2C_2\|\mu_{pp}\|_\tv \ge R$, it is not
possible to meet the condition $\|\omega_0\|_{X_T} \le R/2$ by an 
appropriate choice of $T$, so the argument breaks down and does not 
give any information on the existence of solutions to \eqref{omin}. 

\subsubsection{End of the Proof of Theorem~\ref{thm1}} 
The solution $\omega \in X_T$ of \eqref{omin} constructed by the 
fixed point argument automatically satisfies $\omega \in C^0((0,T],
L^1(\R^2))$. Indeed, arguing as in \eqref{Fest}, we find
\begin{equation}\label{omL1}
  \|\omega(t) - \omega_0(t)\|_{L^1} \,\le\, \int_0^t \frac{C}{
  (t-s)^{\frac12}}\,\|\omega(s)\|_{L^{4/3}}^2\d s \,\le\, 
  C \|\omega\|_{X_T}^2~,
\end{equation}
for $t \in (0,T]$, where $\omega_0(t) = e^{t\Delta}\mu$. Since
$\|\omega_0(t)\|_{L^1} \le \|\mu\|_\tv$ by \eqref{heat1}, we conclude
that $\omega(t) \in L^1(\R^2)$ for all $t \in (0,T]$, and the
continuity with respect to time is again easy to verify. In fact, one
can even show that
\begin{equation}\label{L1approx}
  \lim_{t \to 0}\|\omega(t) - \omega_0(t)\|_{L^1} \,=\, 0~.
\end{equation}
Since $\omega_0(t) \weakto \mu$ as $t \to 0$, this implies that
$\omega(t) \weakto \mu$ as $t \to 0$. Moreover, as $\|\omega(t)\|_{L^1}$ 
is a nonincreasing function of $t$, we deduce from \eqref{L1approx}
that $\|\omega(t)\|_{L^1} \le \|\mu\|_\tv$ for all $t \in (0,T]$.

To prove \eqref{L1approx}, we denote
\[
  \delta \,:=\, \limsup_{t \to 0}\,t^{1/4}\|\omega(t) - \omega_0(t)
  \|_{L^{4/3}} \,\equiv\, \limsup_{T \to 0}\|\omega - \omega_0\|_{X_T}~.
\]
Since $\omega - \omega_0 = (F\omega_0 - F\omega) - F\omega_0$  
and $\|\omega\|_{X_T} + \|\omega_0\|_{X_T} \le 2R$, it follows 
from \eqref{LipF} that $\delta \le (2C_3 R)\delta + \ell_{4/3}(\mu)$, 
where 
\[
  \ell_p(\mu) \,=\, \limsup_{t \to 0}t^{1-\frac1p}
  \|(F\omega_0)(t)\|_{L^p}~, \qquad 1 \le p \le \infty~. 
\]
Now, a direct calculation, which is postponed to Lemma~\ref{ellem}
below, shows that $\ell_p(\mu) = 0$ for any $\mu \in \cM(\R^2)$ and
any $p \in [1,\infty]$. This implies that $\delta = 0$ because 
$2C_3 R < 1$. Using again the identity $\omega - \omega_0 = (F\omega_0 
- F\omega) - F\omega_0$ and arguing as in \eqref{omL1}, we conclude 
that
\[
  \limsup_{t \to 0}\|\omega(t) - \omega_0(t)\|_{L^1} \,\le\, 
  C R \delta + \ell_1(\mu) \,=\, 0~.
\]

To finish the proof of Theorem~\ref{thm1}, it remains to show that all
solutions of \eqref{omeq} are global in time.  We first observe that,
if $\mu \in \cM(\R^2) \cap L^{4/3}(\R^2)$, then
$\|e^{t\Delta}\mu\|_{L^{4/3}} \le \|\mu\|_{L^{4/3}}$ for all $t > 0$,
hence $\|\omega_0\|_{X_T} \le T^{1/4}\|\mu\|_{L^{4/3}}$. Since the
local existence time $T > 0$ is determined by the condition
$\|\omega_0\|_{X_T} \le R/2$, we see that {\em an upper bound} on the
$L^{4/3}$-norm of the initial data provides a {\em lower bound} on the
local existence time $T$. Assume now that $\mu \in \cM(\R^2)$
satisfies $2C_2 \|\mu_{pp}\|_\tv < R$, and let $\omega \in X_T$ be the
local solution of \eqref{omeq} constructed by the fixed point
argument. If $T_1 > 0$ satisfies $T_1^{1/4}\|\omega(T) \|_{L^{4/3}}
\le R/2$, the same argument provides a solution $\tilde \omega$ of
\eqref{omeq} on the time interval $[T,T+T_1]$ with $\tilde \omega(T) =
\omega(T)$, and arguing as in \cite{BA} one can verify that $\tilde
\omega \in C^0([T,T+T_1],L^1(\R^2)\cap L^{4/3}(\R^2))$. Thus, if we
glue together $\omega$ and $\tilde \omega$, we obtain a local solution
$\omega \in C^0((0,T+T_1],L^1(\R^2)\cap L^{4/3}(\R^2))$ of
\eqref{omeq}. Iterating this procedure, and using the fact that
$\|\omega(t)\|_{L^{4/3}}$ is a nonincreasing function of time, see
\eqref{Lpdecay}, we can construct a (unique) local solution on the
time interval $(0,T+kT_1]$ for any $k \in \N$. This means that the
solution of \eqref{omeq} is global if $2C_2 \|\mu_{pp}\|_\tv < R$.

The proof of Theorem~\ref{thm1} is now complete except for two minor
points. First, once we know that the solution $\omega$ of \eqref{omeq}
lies in $L^{4/3}(\R^2)$ for positive times, a standard bootstrap
argument shows that $\omega(t) \in L^\infty(\R^2)$ for $t > 0$, 
as asserted in Theorem~\ref{thm1}. More importantly, while the fixed 
point argument only shows that the solution is unique in a ball of 
the space $X_T$, one can in fact prove uniqueness in the whole space 
$C^0_b((0,T],L^1(\R^2)) \cap C^0((0,T],L^\infty(\R^2))$, for any $T > 0$. 
Establishing this improved uniqueness property requires additional 
arguments which will be presented in Section~\ref{ss3.1} below.

\medskip
Theorem~\ref{thm1} shows that the vorticity equation \eqref{omeq}
is globally well-posed if the atomic part of the initial data satisfies
$\|\mu_{pp}\|_\tv \le C_0 \nu$ for some $C_0 > 0$. As was already 
mentioned, this smallness condition can be completely relaxed, 
and the optimal result is\: 

\begin{theorem}\label{thm2} {\rm \cite{GMO,GG}}\\
For any $\mu\in \cM(\R^2)$, the vorticity equation \eqref{omeq} 
has a unique global solution
\[
  \omega \in C^0((0,\infty),L^1(\R^2) \cap L^\infty(\R^2))
\]
such that $\|\omega(\cdot,t)\|_{L^1} \le \|\mu\|_\tv$ for all 
$t > 0$, and $\omega(\cdot,t) \weakto \mu$ as $t \to 0$.
\end{theorem}

The existence part in Theorem~\ref{thm2} can be established using a
relatively standard approximation procedure, see \cite{GMO}, but
proving uniqueness without any restriction on the initial data
requires a careful treatment of the large Dirac masses in the initial
measure. We refer to Section~\ref{s3} below for a sketch of the
uniqueness proof, and to \cite{GG} for full details.

\subsubsection{Short Time Self-Interaction}

To conclude this section, we state and prove an auxiliary 
result which was used in the proof of Theorem~\ref{thm1}. 

\begin{lemma}\label{ellem} Let $\mu \in \cM(\R^2)$. For any 
$t > 0$, let $\omega_0(t) = e^{t\Delta}\mu$ and let $u_0(t) 
= K*\omega_0(t)$ be the corresponding velocity field, given by the 
Biot-Savart law \eqref{BS}. Then
\begin{equation}\label{ellpzero}
  \ell_p(\mu) \,:=\, \limsup_{t \to 0}\,t^{1-\frac1p}
  \|(F\omega_0)(t)\|_{L^p} \,=\, 0~, \qquad 1 \le p \le \infty~,
\end{equation}
where 
\[
  (F\omega_0)(t) \,=\, \int_0^t e^{(t-s)\Delta} \,u_0(s)\cdot
  \nabla\omega_0(s)\d s~, \qquad t > 0~. 
\]
\end{lemma}

\noindent{\bf Proof.}
We shall prove \eqref{ellpzero} for $p \in 
[1,2)$. The other values of $p$ can be treated in a similar way. 
Arguing as in \eqref{Fest}, we find
\[
  t^{1-\frac1p}\|(F\omega_0)(t)\|_{L^p} \,\le\, 
  t^{1-\frac1p}\int_0^t \frac{C}{(t-s)^{\frac32 - \frac1p}}
  \,\|\omega_0(s)\|_{L^{4/3}}^2\d s~, \quad t > 0~.
\]
It follows that $\ell_p(\mu) \le C_p (L_{4/3}(\mu))^2$, where $C_p$ 
is a constant depending only on $p$ and $L_p(\mu)$ is defined 
in \eqref{heat2}. This already shows that $\ell_p(\mu) = 0$ 
if $\mu_{pp} = 0$. More generally, the argument above implies that 
\begin{equation}\label{ellpdiff}
  |\ell_p(\mu_1 + \mu_2) - \ell_p(\mu_1) - \ell_p(\mu_2)| \,\le\, 
  2C_p L_{4/3}(\mu_1)L_{4/3}(\mu_2)~,
\end{equation}
for any $\mu_1, \mu_2 \in \cM(\R^2)$. Taking $\mu_1 = \mu - \mu_{pp}$
and $\mu_2 = \mu_{pp}$, we conclude that $\ell_p(\mu) =  \ell_p(\mu_{pp})$.   
This means that the non-atomic part of the measure $\mu$ does not 
contribute to the limit $\ell_p(\mu)$. It is thus sufficient to 
prove \eqref{ellpzero} in the case where $\mu$ is purely atomic, 
and in view of \eqref{ellpdiff} we can even suppose that $\mu$ is
a finite linear combination of Dirac masses\: $\mu = \sum_{i=1}^N
\alpha_i \delta_{x_i}$. In that case, 
\[
  \omega_0(x,t) \,=\, \sum_{i=1}^N \frac{\alpha_i}{t}\,G\Bigl(
  \frac{x-x_i}{\sqrt{t}}\Bigr)~, \qquad
  u_0(x,t) \,=\, \sum_{i=1}^N \frac{\alpha_i}{\sqrt{t}}\,v^G\Bigl(
  \frac{x-x_i}{\sqrt{t}}\Bigr)~,
\]
where $G$ and $v^G$ are defined in \eqref{Gdef} below. Let 
\[
  H(x,t) \,=\, u_0(x,t) \cdot \nabla \omega_0(x,t) \,=\, 
  \sum_{i\neq j}\frac{\alpha_i\alpha_j}{t^2}\,v^G\Bigl(
  \frac{x-x_i}{\sqrt{t}}\Bigr) \cdot \nabla G\Bigl(
  \frac{x-x_j}{\sqrt{t}}\Bigr)~.
\]
Here the sum runs over indices $i \neq j$ because the contributions of
the self-interaction terms $i = j$ are identically zero. Using 
this observation, it is rather straightforward to verify that
\[
  \limsup_{t \to 0} t^{\frac32 - \frac1p} \|H(t)\|_{L^p} \,<\, \infty~,
  \qquad 1 \le p \le \infty~.
\]
As $p < \infty$, it follows that
\[
  t^{1-\frac1p}\|(F\omega_0)(t)\|_{L^p} \,\le\, t^{1-\frac1p} \int_0^t
  \frac{C}{(t-s)^{1-\frac1p}}\,\|H(s)\|_{L^1}\dd s 
  \,\xrightarrow[t \to 0]{}\, 0~.
\]
This concludes the proof. 
\QED

\section{Self-Similar Variables and Long-Time Behavior}
\label{s3}

Explicit examples of two-dimensional viscous flows are easily 
constructed if one assumes radial symmetry. Indeed, if the 
vorticity $\omega(x,t)$ is radially symmetric, the velocity 
field $u = K * \omega$ given by the Biot-Savart law \eqref{BS}
is {\em azimuthal}, in the sense that $x \cdot u(x,t) \equiv 0$. 
As a consequence, the nonlinearity $u \cdot \nabla\omega$ in \eqref{omeq}
vanishes identically, and \eqref{omeq} therefore reduces to the 
linear heat equation $\partial_t \omega = \nu\Delta\omega$ 
(which preserves radial symmetry). 

As an example, consider the particular case where the initial 
vorticity is a point vortex of circulation $\gamma \in \R$ located
at the origin\: $\mu = \gamma \delta_0$. The unique solution of 
\eqref{omeq} given by  Theorem~\ref{thm2} is the {\em Lamb-Oseen 
vortex}\:
\begin{equation}\label{Osdef}
  \omega(x,t) \,=\, \frac{\gamma}{\nu t}\,G\Bigl(\frac{x}{\sqrt{\nu t}}
  \Bigr)~,\qquad u(x,t) \,=\, \frac{\gamma}{\sqrt{\nu t}}\,
  v^G\Bigl(\frac{x}{\sqrt{\nu t}}\Bigr)~,
\end{equation}
where the vorticity and velocity profiles are given by
\begin{equation}\label{Gdef}
  G(\xi) \,=\, \frac{1}{4\pi}\,e^{-|\xi|^2/4}~, \qquad 
  v^G(\xi) \,=\, \frac{1}{2\pi}\,\frac{\xi^\perp}{|\xi|^2}
  \Bigl(1 - e^{-|\xi|^2/4}\Bigr)~, \qquad \xi \in \R^2~.
\end{equation}
Note that $\int_{\R^2} G(\xi)\d \xi = 1$, so that $\gamma =
\int_{\R^2} \omega(x,t)\d x$ is the total circulation of the vortex.
Oseen's vortex is thus a {\em self-similar solution} of \eqref{omeq}
with Gaussian vorticity profile $G$. The velocity profile $v^G$ is
azimuthal, vanishes at the origin, and satisfies $|v^G(\xi)| \sim
(2\pi|\xi|)^{-1}$ as $|\xi| \to \infty$. In particular, $v^G \notin
L^2(\R^2)$, in agreement with Lemma~\ref{kinerg}. Oseen vortices are
therefore infinite energy solutions of the two-dimensional
Navier-Stokes equation.

Oseen's vortex plays an important role in the dynamics of the 
vorticity equation \eqref{omeq}. As was already observed, it 
can be considered as the ``fundamental solution'' of \eqref{omeq}; 
i.e., the solution with a single Dirac mass as initial data. 
In addition, the following result shows that it describes 
the long-time asymptotics of all integrable solutions of 
\eqref{omeq}\:

\begin{theorem}\label{thm3} {\rm \cite{GK,GW2}}\\
For all initial data $\mu \in \cM(\R^2)$, the solution $\omega(x,t)$ 
of \eqref{omeq} given by Theorem~\ref{thm2} satisfies 
\begin{equation}\label{Osconv}
   \lim_{t \to \infty} \Big\|\omega(x,t) - \frac{\gamma}{\nu t}
   G\Bigl(\frac{x}{\sqrt{\nu t}}\Bigr)\Big\|_{L^1} \,=\, 0~,
   \qquad \hbox{where} \quad  \gamma \,=\, \int_{\R^2} \d\mu~.
\end{equation}
\end{theorem}

It follows immediately from \eqref{Osconv} that Oseen vortices are the
only self-similar solutions of the Navier-Stokes equation in $\R^2$
for which the vorticity profile is integrable.  As an aside, we
mention that, following the approach of Cannone and Planchon
\cite{CP}, one can construct an infinite-dimensional family of small
self-similar solutions of \eqref{NS} for which $u \in
L^{2,\infty}(\R^2)$ but $\omega \notin L^1(\R^2)$. Assuming, as in
Theorem~\ref{thm3}, that the initial measure is finite allows to
eliminate all these ``artificial'' solutions, which are produced by
fat tails in the initial data and not by the intrinsinc dynamics
of Eq.~\eqref{omeq}. Also, Theorem~\ref{thm3} strongly suggests that
Oseen vortices are {\em stable} solutions of \eqref{omeq} for all
values of the circulation $\gamma$. This important question will be
discussed in Section~\ref{s4} below.

It is worth mentioning that the convergence result \eqref{Osconv} is
not constructive, and does not provide any estimate of the time needed
for the solution of \eqref{omeq} to approach Oseen's vortex.  Under
strong localization assumptions, such estimates were obtained in
\cite{GR}, but these results are certainly not optimal. We also refer
to \cite{CPR} for an interesting attempt toward a more precise
description of the intermediate asymptotics, using ideas from
statistical mechanics.

In the rest of this section, we give a sketch of the proof of
Theorem~\ref{thm3} and of the uniqueness part of Theorem~\ref{thm2}.
We first recall classical estimates for solutions of
convection-diffusion equations, which are due to Osada \cite{Os}.  We
then introduce self-similar variables, which allow to compactify the
solutions of \eqref{omeq} in $L^1(\R^2)$, and to transform the
self-similar Oseen vortices into equilibria. Finally, we use Lyapunov
functions to determine the $\omega$-limit sets of solutions of
\eqref{omeq}, and we establish a ``Liouville theorem'' which implies
both the convergence result \eqref{Osconv} and the uniqueness of the
solution of \eqref{omeq} with a single Dirac mass as initial
vorticity.

\subsection{Estimates for Convection-Diffusion Equations}
\label{ss3.1}

Given $\nu > 0$, we consider a linear convection-diffusion equation 
of the form
\begin{equation}\label{Ueq}
  \partial_t \omega(x,t) + U(x,t) \cdot \nabla \omega(x,t) 
  \,=\, \nu\Delta \omega(x,t)~,
\end{equation}
where $x \in \R^2$, $t \in (0,T)$, and $U : \R^2 \times 
(0,T) \to \R^2$ is a given divergence-free vector field. 
We assume that the curl $\Omega = \partial_1 U_2 - \partial_2 U_1$ 
of the advection field $U$ satisfies $\Omega \in C^0((0,T),L^1(\R^2))$ 
and $\|\Omega(\cdot,t)\|_{L^1(\R^2)} \le K_0 \nu$ for $t \in (0,T)$, where
$K_0$ is a positive constant. Then, in view of \cite{Os}, any solution 
of \eqref{Ueq} can be represented as
\begin{equation}\label{omrepr}
  \omega(x,t) \,=\, \int_{\R^2} \Gamma^\nu_U(x,t;y,s)\omega(y,s)\d y~,
  \quad x \in \R^2~, \quad 0 < s < t < T~,
\end{equation}
where $\Gamma^\nu_U$ is the {\em fundamental solution} of the 
convection-diffusion equation \eqref{Ueq}. The following
properties of $\Gamma^\nu_U$ will be useful:

\medskip\noindent{$\bullet$} There exist $\beta > 0$ and 
$K_1 > 0$ (depending only on $K_0$) such that 
\begin{equation}\label{gamm1}
  0 \,<\, \Gamma^\nu_U(x,t;y,s) \,\le\, \frac{K_1}{\nu(t-s)}
  \,\exp\Bigl(-\beta\frac{|x-y|^2}{4\nu(t-s)}\Bigr)~,
\end{equation}
for $x,y \in \R^2$ and $0 < s < t < T$, see \cite{Os}. In fact, a
more precise result due to Carlen and Loss \cite{CL} shows that, if
$(\nu t)^{1/2} \|U(\cdot,t)\|_{L^\infty} \le K_0 \nu$ for all $t \in
(0,T)$, then \eqref{gamm1} holds for any $\beta \in (0,1)$, with a
constant $K_1$ depending on $K_0$ and $\beta$. It is also possible
to establish a Gaussian lower bound on $\Gamma^\nu_U$, see \cite{Os}.

\medskip\noindent{$\bullet$} There exists $\gamma \in (0,1)$ 
(depending only on $K_0$) and, for any $\delta > 0$, there exists
$K_2 > 0$ (depending only on $K_0$ and $\delta$) such that 
\begin{align}\nonumber
  |\Gamma^\nu_U(x,t;y,s) &- \Gamma^\nu_U(x',t';y',s')| \\ \label{gamm2}
  \,&\le\, K_2\Bigl(|x{-}x'|^\gamma + |\nu(t{-}t')|^{\gamma/2} + 
  |y{-}y'|^\gamma + |\nu(s{-}s')|^{\gamma/2}\Bigr)~,
\end{align}
whenever $\nu(t-s) \ge \delta$ and $\nu(t'-s') \ge \delta$, see~\cite{Os}.

\medskip\noindent{$\bullet$} For $0 < s < t < T$ and $x,y \in
\R^2$,  
\begin{equation}\label{gamm3a}
  \int_{\R^2} \Gamma^\nu_U(x,t;y,s)\d x \,=\, 1~, \quad 
  \int_{\R^2} \Gamma^\nu_U(x,t;y,s)\d y \,=\, 1~. 
\end{equation}
For $0 < s < r < t < T$ and $x,y \in \R^2$,  
\begin{equation}\label{gamm3b}
  \Gamma^\nu_U(x,t;y,s) \,=\, \int_{\R^2} \Gamma^\nu_U(x,t;z,r) 
  \Gamma^\nu_U(z,r;y,s)\d z~.  
\end{equation}

\begin{remark}\label{uptozero}
If $x,y \in \R^2$ and $t > 0$, it follows from~\eqref{gamm2} that 
the function $s \mapsto \Gamma^\nu_U(x,t;y,s)$ can be continuously extended 
up to $s = 0$, and that this extension (still denoted by $\Gamma^\nu_U$) 
satisfies properties~\eqref{gamm1} to~\eqref{gamm3b} with $s = 0$.
\end{remark}

As an application of these results, we give a new formulation of the
uniqueness claims in Theorems~\ref{thm1} and \ref{thm2}.  Assume that
$\omega \in C^0((0,T),L^1(\R^2) \cap L^\infty(\R^2))$ is a solution of
\eqref{omeq} which is uniformly bounded in $L^1(\R^2)$\: there exists
$K_0 > 0$ such that $\|\omega(t)\|_{L^1} \le K_0\nu$ for $t \in
(0,T)$.  Let $\phi \in C_0^\infty(\R^2)$ be a test function. Using
\eqref{omeq} and integrating by parts, we easily find
\[
  \frac{\dd}{\dd t} \int_{\R^2} \omega(x,t)\phi(x)\d x \,=\, 
  \int_{\R^2} \omega(x,t) u(x,t)\cdot \nabla \phi(x)\d x
  + \nu\int_{\R^2} \omega(x,t) \Delta \phi(x)\d x~.  
\]
If we express $u$ in terms of $\omega$ using the Biot-Savart 
law \eqref{BS}, we can rewrite the first integral on the right-hand 
side as
\[
  \frac{1}{4\pi} 
  \int_{\R^2 \times \R^2}\omega(x,t)\omega(y,t) \frac{(x-y)^\perp}{
  |x-y|^2}\cdot\Bigl(\nabla\phi(x) - \nabla\phi(y)\Bigr)\d x \d y~, 
\]
see e.g. \cite{De}. Since $\omega$ is uniformly bounded in $L^1(\R^2)$, 
we conclude that 
\[
  \Bigl|\frac{\dd}{\dd t} \int_{\R^2} \omega(x,t)\phi(x)\d x \Bigr|
  \,\le\, C\nu^2 \sum_{|\alpha| = 2}\|\partial^\alpha \phi\|_{L^\infty}~,
  \qquad t \in (0,T)~,
\]
for some constant $C > 0$ depending only on $K_0$. This implies that
$\omega(t)$ has a limit in $D'(\R^2)$ as $t \to 0$, which we denote by
$\mu$. As $\|\omega(t)\|_{L^1} \le K_0\nu$, it follows that $\mu \in
\cM(\R^2)$ with $\|\mu\|_\tv \le K_0\nu$, and that $\omega(t) \weakto
\mu$ as $t \to 0$.

On the other hand, since $\omega$ solves \eqref{Ueq} with 
$U(x,t) = u(x,t)$, and since the curl of the advection field 
is uniformly bounded in $L^1(\R^2)$, we have the representation 
\eqref{omrepr}, where the fundamental solution $\Gamma^\nu_u(x,t;y,s)$ 
satisfies \eqref{gamm1} to \eqref{gamm3b}. In particular, using 
Remark~\ref{uptozero}, we have for all $x \in \R^2$ and all 
$t \in (0,T)$, 
\begin{align*}
  \omega(x,t) \,&=\,  \int_{\R^2} \Gamma^\nu_u(x,t;y,0)\omega(y,s)\d y \\
  \,&+\, \int_{\R^2} (\Gamma^\nu_u(x,t;y,s) - \Gamma^\nu_u(x,t;y,0))
  \omega(y,s)\d y~, \quad 0 < s < t~.
\end{align*}
In view of \eqref{gamm2}, the second integral on the right-hand side
converges to zero as $s \to 0$. On the other hand, since $y \mapsto 
\Gamma^\nu_u(x,t;y,0)$ is continuous and vanishes at infinity, and 
since $\omega(s) \weakto \mu$ as $s \to 0$, we can take the
limit $s \to 0$ in the first integral and we obtain the useful formula
\begin{equation}\label{representation}
  \omega(x,t) \,=\, \int_{\R^2} \Gamma^\nu_u(x,t;y,0)\d \mu_y~, 
  \qquad x \in \R^2~, \quad 0 < t < T~,
\end{equation}
which shows that any solution $\omega(x,t)$ of \eqref{omeq}
which is uniformly bounded in $L^1(\R^2)$ can be represented
in terms of its trace $\mu$ at $t = 0$, using the corresponding 
fundamental solution $\Gamma^\nu_u$.

In particular, since $\Gamma^\nu_u(x,t;y,0)$ is positive and satisfies 
\eqref{gamm3a}, it follows from \eqref{representation} that 
$\|\omega(t)\|_{L^1} \le \|\mu\|_\tv$ for all $t \in (0,T)$. 
Thus we can assume that $K_0 = \|\mu\|_\tv$ without loss of generality. 
Moreover, using the upper bound \eqref{gamm1} we deduce from 
\eqref{representation} that
\begin{equation}\label{omupp}
  \sup_{0 < t < T} (\nu t)^{1-\frac1p} \|\omega(t)\|_{L^p} \,\le\, 
  C\|\mu\|_\tv~, \qquad 1 \le p \le \infty~,
\end{equation}
where $C$ is a positive constant depending only on $K_0$. In fact, 
a more careful argument shows that the constant $C > 0$ in \eqref{omupp} 
is independent of $K_0$, see \cite{BA,Ka3,CL}. Similarly, repeating 
the proof of Lemma~\ref{heatlem}, we find 
\begin{equation}\label{omupp2}
  \limsup_{t \to 0} (\nu t)^{1-\frac1p} \|\omega(t)\|_{L^p} \,\le\, 
  C \|\mu_{pp}\|_\tv~, \qquad 1 < p \le \infty~. 
\end{equation}
Taking $p = 4/3$, we deduce from \eqref{omupp}, \eqref{omupp2}
that the solution $\omega$ of \eqref{omeq} lies in the space $X_T$ 
defined in \eqref{XTdef}, and that $\|\omega\|_{X_T}$ is small if  
$\|\mu\|_\tv$ is small, or if $\|\mu_{pp}\|_\tv$ and $T$ are small. 
In both cases, $\omega$ necessarily coincides with the solution 
of \eqref{omeq} constructed by the fixed point argument. 

\subsection{Self-Similar Variables}
\label{ss3.2}

Let $\omega \in C^0((0,\infty),L^1(\R^2) \cap L^\infty(\R^2))$ be 
a solution of \eqref{omeq}. Given $x_0 \in \R^2$, $t_0 > 0$, and 
$T > 0$, we introduce the {\em self-similar variables}
\begin{equation}\label{xitaudef}
   \xi \,=\, \frac{x-x_0}{\sqrt{\nu(t+t_0)}}~, \qquad 
   \tau \,=\, \log\Bigl(\frac{t+t_0}{T}\Bigr)~,
\end{equation}
and we transform the vorticity and velocity fields as follows\:
\begin{equation}\label{wvdef}
  \begin{aligned}
  \omega(x,t) \,&=\, \frac{1}{t+t_0} ~w\left(\frac{x-x_0}{\sqrt{\nu
  (t+t_0)}}\,,\,\log\Bigl(\frac{t+t_0}{T}\Bigr)\right)~, \\
  u(x,t) \,&=\, \sqrt{\frac{\nu}{t+t_0}} ~v \left(\frac{x-x_0}{ 
  \sqrt{\nu(t+t_0)}}\,,\,\log\Bigl(\frac{t+t_0}{T}\Bigr)\right)~.
  \end{aligned}
\end{equation}
Then the rescaled vorticity $w(\xi,\tau)$ satisfies the evolution 
equation
\begin{equation}\label{weq}
  \partial_\tau w + v\cdot\nabla_\xi\,w \,=\, \Delta_\xi\,w + \frac12
  \,\xi\cdot\nabla_\xi\,w + w~,
\end{equation}
and the rescaled velocity $v$ is again expressed in terms of $w$ via
the Biot-Savart law \eqref{BS}\: $v = K * w$. The initial data of 
\eqref{omeq} and \eqref{weq} are related via
\[
  w(\xi,\tau_0) \,=\, t_0\,\omega(x_0 + \xi\sqrt{\nu t_0}, 0)~, 
  \qquad \xi \in \R^2~,
\]
where $\tau_0 = \log(t_0/T)$. Note that $w(\xi,\tau)$ and $v(\xi,\tau)$ 
are now dimensionless quantities, as are the new space and time variables 
$\xi$, $\tau$. If $w(\xi,\tau) = \alpha G(\xi)$ for some $\alpha \in \R$,
then by \eqref{wvdef} $\omega(x,t)$ is Oseen's vortex with circulation
$\gamma = \alpha \nu$ located at point $x_0 \in \R^2$ and originating
from time $-t_0$. Thus Oseen vortices are now {\em equilibria} of 
the rescaled system~\eqref{weq}. In what follows, we shall use 
the change of variables \eqref{wvdef} with $x_0 = 0$ and $t_0 = T$,
so that $\tau_0 = 0$, but other choices of $x_0$ or $t_0$ are sometimes 
more appropriate, see e.g. \cite{GR}. 

The main effect of the transformation \eqref{wvdef} is to replace
the diffusion operator $\nu\Delta$ on the right-hand side of 
\eqref{omeq} by the more complicated operator
\begin{equation}\label{LLdef}
  \cL \,=\, \Delta  + \frac12\,\xi\cdot\nabla + 1~,
\end{equation}
which appears in \eqref{weq}. As we shall see in Section~\ref{s4}, 
the operator $\cL$ has better spectral properties than the Laplacian, 
when acting on appropriate function spaces. Moreover, the semigroup
$e^{\tau\cL}$ generated by $\cL$ is {\em asymptotically confining}, 
a property that does not hold for the heat semigroup $e^{\nu t \Delta}$.
In fact, using \eqref{heat} and \eqref{wvdef}, it is straightforward to 
derive the explicit formula
\begin{equation}\label{expcL}
  (e^{\tau\cL}w_0)(\xi) \,=\, \frac{1}{4\pi a(\tau)}\int_{\R^2}
  \exp\Bigl(-\frac{|\xi-\eta\,e^{-\tau/2}|^2}{4a(\tau)}\Bigr) 
  w_0(\eta)\d\eta~,
\end{equation}
for all $\xi \in \R^2$ and all $\tau > 0$, where $a(\tau) = 1 - e^{-\tau}$. 
If $w_0 \in L^1(\R^2)$, then Lebesgue's dominated convergence theorem
shows that $(e^{\tau\cL}w_0)(\xi)$ converges to $\alpha G(\xi)$ as 
$\tau \to \infty$, where $\alpha = \int_{\R^2}w_0\d\xi$ and $G$ is 
given by \eqref{Gdef}. Thus the solution of the linear equation 
$\partial_\tau w = \cL w$ is asymptotically confined like a Gaussian 
function provided the initial data are integrable. 

A similar property holds for the nonlinear equation \eqref{weq}, 
and allows to show that the trajectories of this system are 
compact. 

\begin{lemma}\label{complem}
For any $w_0 \in L^1(\R^2)$, the solution $\{w(\tau)\}_{\tau\ge0}$ of 
\eqref{weq} with initial data $w_0$ is relatively compact in 
$L^1(\R^2)$.
\end{lemma}

\noindent{\bf Proof.}
We first observe that the Cauchy problem for
Eq.~\eqref{weq} is globally well-posed in the space $L^1(\R^2)$.
This can be deduced, via the change of variables \eqref{wvdef}, 
from the corresponding statement for the original equation 
\eqref{omeq}, which is a particular case of Theorem~\ref{thm1}. 
This argument also shows that the $L^1$ norm of the solution
$w(\tau)$ is nonincreasing with time, and that the circulation 
parameter $\alpha = \int_{\R^2}w(\xi,\tau)\d\xi$ is a conserved 
quantity. Now, using the representation \eqref{representation}
of the solution of \eqref{omeq} in terms of the initial data, 
and the bound \eqref{gamm1} on the corresponding fundamental 
solution, we obtain the following estimate for the solution
of \eqref{weq}\:
\begin{equation}\label{wupp}
  |w(\xi,\tau)| \,\le\, \frac{K_1}{a(\tau)} \int_{\R^2} 
  \exp\Bigl(-\beta\frac{|\xi-\eta e^{-\tau/2}|^2}{4a(\tau)}
    \Bigr)|w_0(\eta)|\d \eta~,
\end{equation}
for $\xi \in \R^2$ and $\tau > 0$, where again $a(\tau) = 1 - e^{-\tau}$. 
Here $\beta \in (0,1)$ and $K_1 > 0$ depend only on $\|w_0\|_{L^1}$. 
Using this bound, it is rather easy to show that
\begin{equation}\label{comp1}
  \sup_{\tau \ge 1}\int_{|\xi| \ge R} |w(\xi,\tau)| \d\xi 
  \,\xrightarrow[R \to \infty]{}\, 0~, 
\end{equation}
see \cite[Lemma~2.5]{GW2}. On the other hand, classical estimates
for the derivatives of solutions of \eqref{omeq} ensure that
\[
  \sup_{\tau \ge 1} \sup_{\xi \in \R^2} |\nabla w(\xi,\tau)| 
  \,<\, \infty~,
\]
see \cite{BA,Ka2}. This together with \eqref{comp1} implies that
\begin{equation}\label{comp2}
  \sup_{\tau \ge 1} \sup_{|\eta| \le \delta}\int_{\R^2} 
  |w(\xi-\eta,\tau) - w(\xi,\tau)| \d\xi \,\xrightarrow[\delta 
  \to 0]{}\, 0~. 
\end{equation}
By a well-known criterion of Riesz \cite{RS}, \eqref{comp1} and 
\eqref{comp2} together imply that the trajectory $\{w(\tau)\}_{\tau
\ge1}$ is relatively compact in $L^1(\R^2)$. Since we also know
that $w \in C^0([0,1],L^1(\R^2))$, the conclusion follows. 
\QED

If $w_0 \in L^1(\R^2)$, it follows from Lemma~\ref{complem} 
that the solution $w(\tau)$ of \eqref{weq} with initial data
$w_0$ converges as $\tau \to \infty$, in the $L^1$ topology, 
to the $\omega$-limit set
\[
  \Omega(w_0) \,=\, \bigcap\limits_{T \ge 0} \overline{\{w(\tau)\,|\, \tau 
  \ge T\}}~,
\]
where the overline denotes here the closure in $L^1(\R^2)$. This is a
nonempty compact and connected subset of $L^1(\R^2)$ which is
positively and negatively invariant under the evolution defined by
\eqref{weq}, in the following sense\: for any $\bar w_0 \in
\Omega(w_0)$, there exists a complete trajectory $\bar w \in
C^0(\R,L^1(\R^2))$ of \eqref{weq} such that $\bar w(\tau) \in
\Omega(w_0)$ for all $\tau \in \R$ and $\bar w(0) = \bar w_0$.  In the
next section, we shall prove that $\Omega(w_0)$ consists of a single
point $\{\alpha G\}$, where $\alpha = \int_{\R^2}w_0 \d\xi$.

\subsection{Lyapunov Functions and Liouville's Theorem}
\label{ss3.3}

Lemma~\ref{complem} shows that positive trajectories of 
\eqref{weq} in $L^1(\R^2)$ are relatively compact. The
same conclusion does not apply to negative trajectories 
in general. In fact, the following proposition, which 
is the main result of this section, implies that the 
equilibria associated to Oseen vortices are the only 
negative trajectories of \eqref{weq} which are relatively
compact in $L^1(\R^2)$. 

\begin{proposition}\label{liouville} {\rm (``Liouville's 
theorem'')} If $\{w(\tau)\}_{\tau\in\R}$ is a complete trajectory 
of \eqref{weq} which is relatively compact in $L^1(\R^2)$, 
then there exists $\alpha \in \R$ such that $w(\tau) =  
\alpha G$ for all $\tau\in \R$.
\end{proposition}

\noindent{\bf Proof.}
Let $\cA$ and $\Omega$ denote the $\alpha$-limit set and the
$\omega$-limit set of the trajectory $\{w(\tau)\}_{\tau\in\R}$ in
$L^1(\R^2)$\:
\[
  \cA \,=\, \bigcap\limits_{T \le 0} \overline{\{w(\tau)\,|\, \tau 
  \le T\}}~, \qquad \Omega \,=\, \bigcap\limits_{T \ge 0} 
  \overline{\{w(\tau)\,|\, \tau \ge T\}}~.
\]
Both $\cA$, $\Omega$ are nonempty compact sets in $L^1(\R^2)$.
Moreover, they attract $w(\tau)$ in the sense that
$\mathrm{dist}_{L^1}(w(\tau),\cA) \to 0$ as $\tau\to -\infty$, while
$\mathrm{dist}_{L^1}(w(\tau),\Omega) \to 0$ as $\tau \to
+\infty$. Following \cite{GW2}, to characterize these limit sets more
precisely, we use a pair of Lyapunov functions for Eq.~\eqref{weq}.

\smallskip\noindent{\bf First step.} Our first Lyapunov 
function is just the $L^1$ norm of the solution. Let $\Phi : 
L^1(\R^2) \to \R_+$ be the continuous function defined by
\[
  \Phi(w) \,=\, \int_{\R^2} |w(\xi)| \d\xi~, \qquad w \in 
  L^1(\R^2)~,
\]
and let 
\[
   \Sigma \,=\, \left\{w \in L^1(\R^2) \,\Big|\, \int_{\R^2}
   |w(\xi)|\d \xi = \Big| \int_{\R^2} w(\xi)\d \xi\Big|\right\}~.
\]
In other words, a function $w \in L^1(\R^2)$ belongs to $\Sigma$ if 
and only if $w(\xi)$ has almost everywhere a constant sign.
We have already mentioned that $\Phi$ is nonincreasing along 
trajectories of \eqref{weq}, as a consequence of \eqref{Lpdecay}. 
Now, using the strong maximum principle, it is straightforward
to verify that $\Phi$ is in fact strictly decreasing, except 
along trajectories which lie in $\Sigma$\: if $\{\bar w(\tau)
\}_{\tau\ge0}$ is a solution of \eqref{weq} such that $\Phi(\bar 
w(\tau)) = \Phi(\bar w(0))$ for all $\tau \ge 0$, then 
$\bar w(0) \in \Sigma$ (hence $\bar w(\tau) \in \Sigma$ for 
all $\tau \ge 0$). 

\smallskip By LaSalle's invariance principle \cite{SY}, the
$\alpha$-limit set of the trajectory $\{w(\tau)\}_{\tau\in\R}$ lies in
the neutral set $\Sigma$, because $\cA$ is positively invariant under
the evolution defined by \eqref{weq} and entirely contained in a level
set of the Lyapunov function $\Phi$. The same conclusion applies to
the $\omega$-limit set $\Omega$.  In view of the definition of
$\Sigma$, it follows that $\Phi(\bar w) = |\alpha|$ for all $\bar w
\in \cA$ or $\Omega$, where $\alpha = \int_{\R^2}w(\xi, \tau)\d\xi$ is
the circulation of our solution (which is a conserved quantity).  This
in turn implies that $\Phi(w(\tau)) = |\alpha|$ for all $\tau \in \R$,
which is possible only if the whole trajectory
$\{w(\tau)\}_{\tau\in\R}$ lies in the neutral set $\Sigma$.  Thus any
relatively compact complete trajectory of \eqref{weq} has necessarily
a definite sign, which is the sign of the circulation parameter
$\alpha$. If $\alpha = 0$, this means that $w(\xi,\tau) \equiv 0$, in
which case the proof is complete. If $\alpha < 0$, we observe that the
substitution $w(\xi_1,\xi_2,\tau) \mapsto -w(\xi_2,\xi_1,\tau)$ leaves
Eq.~\eqref{weq} unchanged, but reverses the sign of the
circulation. Thus we can assume henceforth that $\alpha > 0$. In that
case $w(\xi,\tau) > 0$ for all $\xi \in \R^2$ and all $\tau \in \R$,
and proceeding as in Section~\ref{ss3.1} one can show that
$w(\xi,\tau)$ is bounded from above and from below by time-independent
Gaussian functions.

\smallskip\noindent{\bf Second step.} As a second Lyapunov
function, we use the {\em relative entropy} of the vorticity 
distribution $w(\tau)$ with respect to the Gaussian $G$ \cite{Vil}.  
If $w : \R^2 \to \R_+$ is a (measurable) positive function with a 
Gaussian upper bound, we define
\begin{equation}\label{Hdef}
  H(w) \,=\, \int_{\R^2} w(\xi) \log \Bigl(\frac{w(\xi)}{G(\xi)}\Bigr) 
  \d\xi~.
\end{equation}
Then $H$ is nonincreasing along trajectories of \eqref{weq},
because 
\begin{equation}\label{H-dot}
  \frac{\d}{\d\tau}H(w(\tau)) \,=\, -\int_{\R^2} w(\xi,\tau) \Big|\nabla
  \log\Bigl(\frac{w(\xi,\tau)}{G(\xi)}\Bigr)\Big|^2 \d\xi \,\le\, 0~.
\end{equation}
To prove \eqref{H-dot}, we compute
\[
  \frac{\d}{\d\tau}H(w(\tau)) \,=\,
  \int_{\R^2} \Bigl(1+\log\frac{w}{G}\Bigr) \partial_\tau 
  w\d\xi = \int_{\R^2} \Bigl(1+\log\frac{w}{G}\Bigr)(\cL w - 
  v\cdot \nabla w) \d\xi~.
\]
Using the identity $\cL w = \div(G \nabla (\frac{w}{G}))$ 
and integrating by parts, we obtain
\begin{align*}
   \int_{\R^2} \Bigl(1+\log\frac{w}{G}\Bigr)(\cL w)\d\xi 
   \,&=\, -\int_{\R^2} \nabla\Bigl(\log\frac{w}{G}\Bigr)\cdot \frac{G}{w}
    \nabla \Bigl(\frac{w}{G}\Bigr) w \d\xi \\
   \,&=\, -\int_{\R^2} w \Big|\nabla\Bigl(\log\frac{w}{G}\Bigr)\Big|^2 
   \d\xi~. 
\end{align*}
On the other hand, using $v\cdot\nabla w = \div (v w)$ and 
integrating by parts, we find
\begin{align*}
  \int_{\R^2} \Bigl(1+\log\frac{w}{G}\Bigr)(v \cdot \nabla w) \d\xi
    \,&=\, \int_{\R^2} (1+\log(4\pi w)) (v \cdot \nabla w) \d\xi
    + \int_{\R^2} \frac{|\xi|^2}{4}(v \cdot \nabla w) \d\xi\ \\
  \,&=\, -\int_{\R^2} v\cdot\nabla w \d\xi - \frac12 \int_{\R^2} 
    (\xi\cdot v) w \d\xi \,=\, 0~.
\end{align*}
Note that the last integral vanishes because $v = K*w$, see
\eqref{secmomid}. This concludes the proof of inequality
\eqref{H-dot}, which shows in addition that $H$ is strictly decreasing
along trajectories of \eqref{weq}, except on the line of equilibria
corresponding to Oseen vortices. More precisely, if $\{\bar
w(\tau)\}_{\tau\ge0}$ is a positive solution of \eqref{weq} (with a
Gaussian upper bound) such that $H(\bar w(\tau)) = H(\bar w(0))$ for
all $\tau \ge 0$, then $\bar w(0) \in \cT$, where $\cT = \{\beta
G\}_{\beta > 0}$ is the half-line of positive equilibria.

\smallskip We now return to the analysis of our complete trajectory
$\{w(\tau)\}_{\tau\in\R}$ of \eqref{weq}. We know from the first step
that $w(\tau) \in \Sigma$ for all $\tau\in \R$, and without loss of
generality we assume that the circulation parameter $\alpha =
\int_{\R^2}w(\xi,\tau)\d\xi$ is positive, which implies that
$w(\xi,\tau)$ is positive and has a Gaussian upper bound. By LaSalle's
invariance principle, both the $\alpha$-limit set $\cA$ and the
$\omega$-limit set $\Omega$ are contained in the line of equilibria
$\cT$, because these sets are positively invariant and contained in a
level set of the Lyapunov function $H$. The conservation of the total
circulation then implies that $\cA = \Omega = \{\alpha G\}$, hence $H
= \alpha \log(\alpha)$ on $\cA$ and $\Omega$. It follows that
$H(w(\tau)) = \alpha\log(\alpha)$ for all $\tau \in \R$, which is
possible only if $w(\tau) \in \cT$ for all $\tau \in \R$. But $\cT$ is
a line of equilibria, hence we must have $w(\tau) = \alpha G$ for all
$\tau \in \R$, as stated in Proposition~\ref{liouville}.
\QED

As a first application of Proposition~\ref{liouville}, we show
that Oseen vortices describe the long-time asymptotics of all 
solutions of \eqref{omeq} in $L^1(\R^2)$.

\medskip\noindent{\bf Proof of Theorem~\ref{thm3}.}
Assume that $\omega \in C^0((0,\infty), L^1(\R^2) \cap
L^\infty(\R^2))$ is a solution of \eqref{omeq} with initial data $\mu
\in \cM(\R^2)$, in the sense of Theorem~\ref{thm2}, and let
$w(\xi,\tau)$ be the rescaled vorticity obtained from $\omega(x,t)$
via the change of variables \eqref{wvdef}, with an arbitrary choice of
the center $x_0 \in \R^2$ and the initial time $t_0 = T > 0$. Then $w
\in C^0((0,\infty),L^1(\R^2) \cap L^\infty(\R^2))$ is a solution of
the rescaled vorticity equation \eqref{weq} such that
\begin{equation}\label{alphagam}
  \alpha = \int_{\R^2}w(\xi,\tau)\d \xi \,=\, \frac{\gamma}{\nu}~,
  \qquad \hbox{where} \quad \gamma \,=\, \int_{\R^2}\omega(x,t)\d x \,=\, 
  \int_{\R^2}\d\mu~.
\end{equation}
By Lemma~\ref{complem}, the trajectory $\{w(\tau)\}_{\tau\ge1}$ 
is relatively compact in $L^1(\R^2)$. Let $\Omega$ denote the 
$\omega$-limit set of this trajectory in $L^1(\R^2)$. Then 
$\Omega$ is a nonempty compact subset of $L^1(\R^2)$ which is 
positively and negatively invariant under the evolution defined
by \eqref{weq}. In other words, $\Omega$ is a collection of 
complete, relatively compact trajectories of \eqref{weq}. 
By Liouville's theorem, all such trajectories are equilibria
of the form $\beta G$, for some $\beta \in \R$. But the 
conservation of the total circulation implies that $\beta = \alpha$, 
so we conclude that $\Omega = \{\alpha G\}$. This exactly means
that $w(\tau) \to \alpha G$ in $L^1(\R^2)$ as $\tau \to \infty$. 
Returning to the original variables, we obtain \eqref{Osconv}. 
\QED

\begin{remark}\label{xtchoice}
Under the assumptions of Theorem~\ref{thm3}, the above proof
shows that
\begin{equation}\label{Osconv2}
   \lim_{t \to \infty} \Big\|\omega(x,t) - \frac{\gamma}{\nu(t+t_0)}
   \,G\Bigl(\frac{x-x_0}{\sqrt{\nu (t+t_0)}}\Bigr)\Big\|_{L^1} \,=\, 0~,
\end{equation}
for any $x_0 \in \R^2$ and any $t_0 > 0$. This might be surprising
at first sight, but there is no contradiction because, as is 
easily verified, two Oseen vortices with different values of 
$x_0$ or $t_0$ converge to each other in the $L^1$ topology as 
$t \to \infty$. However, the choice of the vortex center and
the initial time becomes important if one wants to determine
the optimal decay rate in \eqref{Osconv2}. This is possible if one
assumes that the vorticity $\omega(x,t)$ decays to zero 
sufficiently fast as $|x| \to \infty$, see \cite{GW2,GR} and 
Section~\ref{s4} below. 
\end{remark}

We next show that Proposition~\ref{liouville} implies the
uniqueness of the solution of the vorticity equation 
\eqref{omeq} with a point vortex as initial data.

\begin{proposition}\label{Dirac} {\rm \cite{GW2,GGL}}\\
Let $\omega \in C^0((0,T),L^1(\R^2)\cap L^\infty(\R^2))$
be a solution of \eqref{omeq} which is bounded in $L^1(\R^2)$ 
and satisfies $\omega \weakto \gamma \delta_0$ as $t \to 0$, for 
some $\gamma \in \R$. Then
\begin{equation}\label{omexp}
   \omega(x,t) \,=\, \frac{\gamma}{\nu t} \,G\Bigl(\frac{x}{
   \sqrt{\nu t}}\Bigr)~, \qquad x \in \R^2~,\quad 0 < t < T~.
\end{equation}
\end{proposition}

\noindent{\bf Proof.}
Applying \eqref{representation} in the particular case where 
$\mu = \gamma \delta_0$, we obtain $\omega(x,t) = \gamma 
\Gamma^\nu_u(x,t;0,0)$, where $\Gamma^\nu_u$ is the fundamental 
solution of \eqref{Ueq} with advection field $u = K*\omega$ given by the 
Biot-Savart law \eqref{BS}. In view of \eqref{gamm1}, we thus have
\begin{equation}\label{ombdd}
  |\omega(x,t)| \,\le\, \frac{|\gamma| K_1}{\nu t}
  \,\exp\Bigl(-\beta\frac{|x|^2}{4\nu t}\Bigr)~,
  \qquad x \in \R^2~, \quad 0 < t < T~.
\end{equation}
for some $K_1 > 0$ and $\beta \in (0,1)$. We now apply the change of
variables \eqref{wvdef} with $x_0 = 0$ and $t_0 = 0$.  Since $\tau =
\log(t/T)$, the rescaled vorticity $w(\xi,\tau)$ is then defined for
all $\xi \in \R^2$ and all $\tau < 0$, so that $w \in
C^0((-\infty,0),L^1(\R^2) \cap L^\infty(\R^2))$ is a negative
trajectory of \eqref{weq}. Moreover, the bound \eqref{ombdd} implies
that $|w(\xi,\tau)| \le (K_1|\gamma|/\nu) e^{-\beta |\xi|^2/4}$, for
all $\xi \in \R^2$ and all $\tau < 0$. As in the proof of
Lemma~\ref{complem}, this bound implies that the trajectory
$\{w(\tau)\}_{\tau < 0}$ is relatively compact in $L^1(\R^2)$.
Applying Proposition~\ref{liouville}, we conclude that $w(\xi,\tau) 
= \alpha G(\xi)$ for some $\alpha \in \R$, and the conservation of 
the total circulation implies that $\alpha = \gamma/\nu$. Thus 
$w(\xi,\tau) = (\gamma/\nu)G(\xi)$ for all $\xi \in \R^2$ and all 
$\tau < 0$, and returning to the original variables we obtain 
\eqref{omexp}.
\QED

It is important to emphasize that there is no restriction on the
circulation $\gamma$ in Proposition~\ref{Dirac}.  Thus the conclusion
cannot be obtained by a standard application of Gronwall's lemma,
which would require that the ratio $|\gamma|/\nu$ be sufficiently
small.  Intuitively, Proposition~\ref{Dirac} and Theorem~\ref{thm1}
together imply that the solution of \eqref{omeq} with arbitrary
initial data $\mu \in \cM(\R^2)$ should be unique, as asserted in
Theorem~\ref{thm2}. Indeed, given any $\epsilon > 0$, we can decompose
$\mu = \mu_1 + \mu_2$, where $\mu_1$ is a finite superposition of
Dirac masses and $\|\mu_2\|_\tv \le \epsilon$. If $\epsilon$ is
sufficiently small, we know from Theorem~\ref{thm1} that \eqref{omeq}
has a unique solution with initial data $\mu_2$, and
Proposition~\ref{Dirac} implies that each atom in $\mu_1$ also
generates a unique solution of \eqref{omeq}. Although the 
the vorticity equation \eqref{omeq} is nonlinear, this decomposition
can be used to prove uniqueness if we can show, as in the 
particular case considered in Lemma~\ref{ellem}, that the 
various components of the solution do not strongly interact for 
small times. These rough ideas can be turned into a rigorous 
proof, and we refer the interested reader to \cite{GG}
for more details. 

\section{Asymptotic Stability of Oseen Vortices}
\label{s4}

As was mentioned in Section~\ref{ss3.2}, the rescaled vorticity
equation 
\begin{equation}\label{weq2}
  \partial_\tau w + v\cdot\nabla_\xi\,w \,=\, \cL w~, \qquad 
  \hbox{where}\quad \cL \,=\, \Delta + \frac12\,\xi\cdot\nabla 
  + 1~,
\end{equation}
has a line of equilibria $\{\alpha G\}_{\alpha \in \R}$ which
correspond, in the original variables, to the family of self-similar
Oseen vortices \eqref{Osdef}. In fact, Proposition~\ref{liouville}
implies that these are the only steady states of \eqref{weq2} in
the space $L^1(\R^2)$, and Theorem~\ref{thm3} even shows that the
family  $\{\alpha G\}_{\alpha \in \R}$ is globally attracting in the sense 
that any solution of \eqref{weq2} in $L^1(\R^2)$ converges to 
$\alpha G$ as $\tau \to \infty$, where $\alpha = \int_{\R^2}w(\xi,\tau)
\d\xi$ is the circulation parameter. This strongly suggests, but does
not rigorously prove, that $\alpha G$ is a stable equilibrium of 
\eqref{weq2} for any $\alpha \in \R$. 

\subsection{Stability with Respect to Localized Perturbations}
\label{ss4.1}

The aim of this section is to study the flow of \eqref{weq2}
in a neighborhood of Oseen's vortex $\alpha G$, and to show that
this steady state is indeed stable for any value of the circulation 
parameter $\alpha \in \R$. Setting $w = \alpha G + \tilde w$, 
$v = \alpha v^G + \tilde v$, we obtain the perturbation equation
\begin{equation}\label{weq3}
   \partial_\tau \tilde w + \tilde  v \cdot \nabla \tilde w \,=\, 
   (\cL - \alpha \Lambda)\tilde w~, 
\end{equation}
where
\begin{equation}\label{Lamdef}
  \Lambda \tilde w \,=\, v^G \cdot\nabla \tilde w + 
  \tilde v\cdot\nabla G~.
\end{equation}  
Here $\tilde v = K * \tilde w$ is the velocity field obtained from
$\tilde w$ via the Biot-Savart law \eqref{BS}. The perturbation
equation \eqref{weq3} is globally well-posed in $L^1(\R^2)$, and it
can be proved that the origin $\tilde w = 0$ is a stable equilibrium,
but in such a large function space it is impossible to obtain a more
precise description of the long-time behavior of the solutions.
Indeed, it follows from the results of \cite{GW1,GW2} that the
spectrum of the linearized operator $\cL - \alpha\Lambda$ in
$L^1(\R^2)$ is exactly the full left-half plane $\{z \in \C\,|\,
\Re(z) \le 0\}$, for any value of $\alpha \in \R$. This means that
there exist perturbations $\tilde w$ which converge to zero at an
arbitrarily slow rate.

For a more precise stability analysis, we have to assume that the
perturbations $\tilde w$ have a faster decay as $|\xi| \to \infty$
than what is strictly needed for integrability. From now on, we write
$w,v$ instead of $\tilde w,\tilde v$, and we consider the perturbation
equation \eqref{weq3} in the weighted space $X =
L^2(\R^2,G^{-1}\d\xi)$ equipped with the scalar product
\begin{equation}\label{Xscalar}
  \langle w_1\,,\,w_2\rangle \,=\, \int_{\R^2} G(\xi)^{-1} \,w_1(\xi) 
  w_2(\xi)\d\xi~,
\end{equation}
and the associated norm $\|w\|_X = \langle w,w\rangle^{1/2}$.
Functions in $X$ are locally in $L^2$ but have a Gaussian decay at
infinity; in particular $X \hookrightarrow L^p(\R^2)$ for all $p \in
[1,2]$.  This function space is convenient because, as we shall see,
the linearized operator $\cL - \alpha\Lambda$ has very nice spectral
properties in $X$. However, if we do not not want to restrict
ourselves to perturbations with Gaussian decay at infinity, it is
possible to use the polynomially weighted space $L^2(m) =
L^2(\R^2,(1+|\xi|^2)^m\d\xi)$, for some $m > 1$ \cite{GW2}.

The choice of the Gaussian space $X$ is justified by the following 
result\:

\begin{lemma}\label{adjoint} {\rm \cite{GW1,GW2,Ma}}\\
i) The operator $\cL$ is {\em selfadjoint} in the space $X = L^2(\R^2,
G^{-1}\d\xi)$, with compact resolvent and purely discrete spectrum\:
\begin{equation}\label{LLspec}
  \sigma(\cL) \,=\, \Bigl\{-\frac{n}{2}\,\Big|\, n = 0,1,2,\dots\Bigr\}~.
\end{equation}
ii) The operator $\Lambda$ is {\em skew-symmetric} in the same space 
$X$\: 
\begin{equation}\label{Lamasym}
  \langle \Lambda w_1\,,\,w_2\rangle + \langle w_1\,,\,\Lambda 
  w_2\rangle \,=\, 0~, \qquad \hbox{for all} \quad w_1,w_2 \in 
  D(\Lambda) \subset X~.
\end{equation}
\end{lemma}

\noindent{\bf Proof.}
To prove i), we conjugate $\cL$ with the square root of the 
weight $G$ to obtain the operator
\begin{equation}\label{Lconj}
  L \,=\, G^{-1/2} \,\cL ~G^{1/2} \,=\, \Delta \,-\, \frac{|\xi|^2}{16}
  \,+\, \frac12~,
\end{equation}
which we have to consider as acting on $L^2(\R^2)$. Clearly $L$ is the
two-dimensional harmonic oscillator, which is known to be self-adjoint
in $L^2(\R^2)$ with compact resolvent and discrete spectrum given by
\eqref{LLspec}. Returning to the original operator, we obtain the
desired conclusions, together with the following characterization of
the domain of $\cL$\:
\[ D(\cL) \,=\, \Bigl\{w \in X \,\Big|\, \Delta w \in X\,,~
   \xi \cdot \nabla w \in X\Big\}~.
\]
The proof of ii) is a direct calculation. By \eqref{Lamdef} we 
have $\Lambda = \Lambda_1 + \Lambda_2$, where $\Lambda_1 w = v^G \cdot 
\nabla w$ and $\Lambda_2 w = v \cdot \nabla G = (K * w) \cdot \nabla G$. 
If $w_1, w_2 \in D(\cL) \subset X$, then
\begin{align*}
  \langle \Lambda_1 w_1,w_2\rangle +  \langle w_1, \Lambda_1 w_2
  \rangle \,&=\, \int_{\R^2} G^{-1}\Bigl(w_2 \,v^G \cdot \nabla w_1
  + w_1\,v^G\cdot \nabla w_2\Bigr)\d\xi \cr
  \,&=\, \int_{\R^2} G^{-1} \,v^G\cdot \nabla(w_1 w_2)\d \xi 
  \,=\, 0~, 
\end{align*}
because $G^{-1}v^G$ is divergence-free. Moreover, since $\nabla G 
= -\frac12\xi G$, we have
\begin{align*}
  &\langle \Lambda_2 w_1,w_2\rangle +  \langle w_1, \Lambda_2 w_2
  \rangle \,=\, -\frac12\int_{\R^2} \Bigl((\xi\cdot v_1) w_2 
  + (\xi\cdot v_2) w_1\Bigr)\d\xi \cr
  \,&=\, -\frac14\int_{\R^2}\int_{\R^2}
  \left\{\xi\cdot \frac{(\xi-\eta)^\perp}{|\xi-\eta|^2} +
  \eta \cdot \frac{(\eta-\xi)^\perp}{|\xi-\eta|^2}\right\}
  w_1(\eta) w_2(\xi)\d\eta\d\xi \,=\, 0~,
\end{align*}
see also \eqref{secmomid}. Thus $ \langle \Lambda w_1,w_2\rangle +
\langle w_1, \Lambda w_2 \rangle = 0$ for all $w_1,w_2 \in D(\cL)
\subset X$. By density, this relation holds for all $w_1, w_2$ in the
domain of $\Lambda$, which can be characterized precisely
\cite{Ma}. One can also show that $\Lambda$ is not only
skew-symmetric, but also skew-adjoint in $X$, see \cite{Ma}.
\QED

It is useful to list the first few eigenfunctions of $\cL$ in 
$X$, and to give explicit formulas for the corresponding spectral 
projections. 

\smallskip\noindent 1) The first eigenvalue $\lambda_0 = 0$ 
is simple, and the associated eigenfunction is of course 
the Oseen vortex profile\: $\cL G = 0$. The corresponding 
spectral projection $P_0 : X \to X$ is given by $P_0 w = 
G \int_{\R^2}w\d\xi$. We denote by $X_0 \subset X$ the kernel
of $P_0$, namely
\begin{equation}\label{X0def}
  X_0 \,=\, \Bigl\{w \in X\,\Big|\, \int_{\R^2}w\d \xi = 0\Bigr\}
  \,=\, \{G\}^\perp~,
\end{equation}
where $\{G\}^\perp$ denotes the set of all $w \in X$ which are 
orthogonal to $G$ in $X$. Due to the conservation of the total 
circulation, this subspace of $X$ is invariant under the evolution
defined by the full equation \eqref{weq3}. 

\smallskip\noindent 2) The second eigenvalue $\lambda_1 = -1/2$ 
has multiplicity two, and the associated eigenfunctions are 
the first order derivatives $\partial_1 G$ and $\partial_2 G$.
The spectral projection $P_1 : X \to X$ satisfies $P_1 w = 
-\partial_1G \int_{\R^2}\xi_1 w\d\xi -\partial_2G \int_{\R^2}\xi_2 
w\d\xi$. Let $X_1 = \ker(P_0) \cap \ker(P_1)$, namely
\[
  X_1 \,=\, \Bigl\{w \in X_0\,\Big|\, \int_{\R^2} \xi_i w\d \xi = 0
  \quad\hbox{for}~i = 1,2\Bigr\} \,=\, \{G; \partial_1 G; \partial_2 G
  \}^\perp~.
\]
Due to the conservation of the first-oder moments of the vorticity, 
see \eqref{firstmom}, the subspace $X_1$ is also invariant under 
the evolution defined by \eqref{weq3}. 

\smallskip\noindent 3) The third eigenvalue $\lambda_2 = -1$ 
has multiplicity three, and the associated eigenfunctions are 
the second order derivatives $\Delta G$, $(\partial_1^2 - 
\partial_2^2)G$, and $\partial_1\partial_2 G$. Let $\tilde P_2$ be 
the spectral projection associated only with the symmetric 
eigenfunction $\Delta G$, i.e. $\tilde P_2 w = \frac14 \Delta G
\int_{\R^2}|\xi|^2 w\d\xi$. We denote $X_2 = \ker(P_0) \cap \ker(P_1)
\cap \ker(\tilde P_2)$, namely
\[
  X_2 \,=\, \Bigl\{w \in X_1\,\Big|\, \int_{\R^2} |\xi|^2 w\d \xi = 0
  \Bigr\} \,=\, \{G; \partial_1 G; \partial_2 G; \Delta G\}^\perp~.
\]
In view of \eqref{secondmom}, the subspace $X_2$ is also invariant under 
the evolution defined by \eqref{weq3}. 

\smallskip More generally, for any $n \in \N$, the $n^{\mathrm{th}}$
eigenvalue $\lambda_n = -n/2$ of $\cL$ has multiplicity $n+1$, and the
corresponding eigenspace is spanned by the Hermite functions of degree
$n$ which can be expressed as homogeneous $n^{\mathrm{th}}$ order
derivatives of the Gaussian profile $G$.

We now consider the spectrum of the linearized operator $\cL -
\alpha\Lambda$ in the (complexified) Hilbert space $X$, for any fixed
$\alpha \in \R$. As is easily verified, the operator $\Lambda$ is a
relatively compact perturbation of $\cL$ in $X$. Indeed, if we
decompose $\Lambda = \Lambda_1 + \Lambda_2$ as in the proof of
Lemma~\ref{adjoint}, we see that $\Lambda_2$ is a compact operator in
$X$, while $\Lambda_1 = v^G\cdot\nabla$ is relatively compact with
respect to $\cL$ because this is a first order differential operator
whose coefficients decay to zero as $|\xi| \to \infty$.  By a
well-known perturbation argument \cite{Ka1}, this implies that $\cL -
\alpha\Lambda$ has itself a compact resolvent in $X$.  In particular,
its spectrum is a sequence of (complex) eigenvalues
$\{\lambda_n(\alpha)\}_{n \in \N}$ with finite multiplicities, which can
accumulate only at infinity. In fact, it is not difficult to verify
that $\Re(\lambda_n(\alpha)) \to -\infty$ as $n \to \infty$, for any
$\alpha \in \R$.
 
In general, the eigenvalues of $\cL - \alpha\Lambda$ depend in a 
nontrivial way on the circulation parameter $\alpha$ and are 
not explicitly known. However, if $w \in X$ is a radially 
symmetric eigenfunction of $\cL$, then obviously $\Lambda w = 0$, 
hence $w$ is also an eigenfunction of $\cL - \alpha\Lambda$ for 
any $\alpha \in \R$, and the corresponding eigenvalue does 
not depend on $\alpha$. In particular, $\lambda_0 = 0$ is 
an eigenvalue of $\cL - \alpha\Lambda$ with eigenfunction $G$,
and the orthogonal subspace $X_0 = \{G\}^\perp$ is invariant 
under the linear evolution generated by $\cL - \alpha\Lambda$. 
Moreover, if we differentiate the identity $v^G \cdot \nabla G = 0$ 
with respect to $\xi_1$ and $\xi_2$, we see that $\Lambda (\partial_i
G) = 0$ for $i = 1,2$. This implies that $\lambda_1 = -1/2$ is 
a double eigenfunction of $\cL - \alpha\Lambda$ for any $\alpha
\in \R$, with eigenfunctions $\partial_1 G$ and $\partial_2 G$. 
In addition, the orthogonal subspace $X_1 = \{G;\partial_1 G;
\partial_2 G\}^\perp$ is invariant under the linear evolution 
generated by $\cL - \alpha\Lambda$. Finally, since $\Lambda 
(\Delta G) = 0$ by symmetry, we see that $\lambda_2 = -1$ is 
an eigenvalue of $\cL - \alpha\Lambda$ for any $\alpha \in \R$,
and the subspace $X_2 = \{G;\partial_1 G;\partial_2 G;\Delta G\}^\perp$
is also invariant under the linearized evolution.  

The following simple but important result shows that the equilibria of 
\eqref{weq2} corresponding to Oseen vortices are {\em spectrally stable} 
with respect to perturbations in the Gaussian space $X$.  

\begin{proposition}\label{specstab} {\rm \cite{GW2}}\\
For any $\alpha \in \R$, the spectrum of the linearized operator
$\cL - \alpha\Lambda$ in the space $X$ satisfies
\begin{equation}\label{sigma0}
  \sigma(\cL-\alpha\Lambda) \,\subset\, \Bigl\{z \in \C
  \,\Big|\, \Re(z) \le 0\Bigr\}~.
\end{equation}
Moreover, 
\begin{equation}\label{sigma1}
  \begin{aligned}
  \sigma(\cL-\alpha\Lambda) \,&\subset\, \Bigl\{z \in \C
  \,\Big|\, \Re(z) \le -\frac12\Bigr\} \quad \hbox{in}\quad X_0~, \cr
  \sigma(\cL-\alpha\Lambda) \,&\subset\, \Bigl\{z \in \C
  \,\Big|\, \Re(z) \le -1\Bigr\} \quad \,\hbox{in}\quad X_1~.
  \end{aligned}
\end{equation}
\end{proposition}

\noindent{\bf Proof.}
We already know that the spectrum of $\cL- \alpha\Lambda$ consists
entirely of eigenvalues. Assume that $(\cL-\alpha\Lambda)w = \lambda
w$ for some $\lambda \in \C$ and some normalized vector $w \in D(\cL)
\subset X$. Then, using Lemma~\ref{adjoint}, we find
\begin{equation}\label{dissip}
  \Re(\lambda) \,=\, \Re \langle (\cL-\alpha\Lambda)w,w\rangle
  \,=\, \langle \cL w,w\rangle \,\le\, 0~,
\end{equation}
which proves \eqref{sigma0}. Moreover $\langle \cL w,w\rangle \,\le\, 
-1/2$ if $w \in X_0$ and $\langle \cL w,w\rangle \,\le\, -1$ if 
$w \in X_1$, and \eqref{sigma1} follows. 
\QED

\begin{remark}\label{X2}
If we restrict ourselves to the invariant subspace $X_2$, it 
is easy to show, as in the proof of Proposition~\ref{specstab}, 
that all eigenvalues of the linearized operator $\cL - \alpha\Lambda$ 
satisfy $\Re(\lambda) < -1$ if $\alpha \neq 0$, but it is difficult
to give an optimal upper bound on the real part of the spectrum. 
\end{remark}

As an immediate consequence, we can show that the Oseen vortices
are {\em linearly stable} with respect to perturbations in $X$\:

\begin{corollary}\label{linstab}
For all $\alpha \in \R$, the semigroup generated by the linearized
operator $\cL - \alpha\Lambda$ satisfies
\begin{equation}\label{linest}
  \|e^{\tau(\cL-\alpha\Lambda)}\|_{Z \to Z} \,\le\, e^{-\mu \tau}~,
  \qquad \hbox{for all}\quad \tau \ge 0~,
\end{equation}
where $\mu = 0$ if $Z = X$, $\mu = 1/2$ if $Z = X_0$, and $\mu = 1$ 
if $Z = X_1$. 
\end{corollary}

\noindent{\bf Proof.}
We know that the self-adjoint operator $\cL$ is 
the generator of an analytic semigroup in $X$, and since $\Lambda$ is a 
relatively compact perturbation of $\cL$ it is clear that $\cL - 
\alpha\Lambda$ is also the generator of an analytic semigroup, for any 
$\alpha \in \R$. Now, estimate \eqref{dissip} means that $\cL - 
\alpha\Lambda$ is $m$-dissipative in $X$ \cite{Ka1}. By the 
Lumer-Philips theorem \cite{Pa}, the associated semigroup satisfies
the bound \eqref{linest} with $Z = X$ and $\mu = 0$. Applying the same
argument to the operator $\cL - \alpha\Lambda + \mu$ restricted to
$X_0$ (if $\mu = 1/2$) or $X_1$ (if $\mu = 1$), we obtain the desired
result in the other cases.
\QED

With some additional work, one can also prove that the Oseen 
vortices are {\em nonlinearly stable} with respect to pertubations
in $X$. Here we can restrict ourselves, without loss of generality, 
to perturbations with zero mean; i.e., we can study Eq.~\eqref{weq3}
in the invariant subspace $X_0 \subset X$ defined by \eqref{X0def}. 
Indeed, adding a perturbation with nonzero mean to the equilibrium 
$\alpha G$ is equivalent to adding a perturbation with zero mean 
to some modified equilibrium $\tilde \alpha G$,  with $\tilde \alpha$ 
close to $\alpha$. The result is\:

\begin{proposition}\label{locstab} {\rm \cite{GW2,GR}}\\
There exists $\epsilon > 0$ such that, for any $\alpha \in \R$, we 
have the following result. If $w_0 \in X_0$ satisfies $\|w_0\|_X \le 
\epsilon$, then the unique solution of \eqref{weq3} with initial data 
$w_0$ satisfies
\begin{equation}\label{wdecay}
  \|w(\tau)\|_X \,\le\, \min(1,2 e^{-\mu \tau})\|w_0\|_X ~, \qquad 
  \hbox{for all} \quad \tau \ge 0~, 
\end{equation}
where $\mu = 1$ if $w_0 \in X_1$ and $\mu = 1/2$ if $w_0 \in 
X_0 \setminus X_1$. 
\end{proposition}

\noindent{\bf Proof.}
Arguing as in Section~\ref{ss2.3}, it is not difficult to verify that
the perturbation equation \eqref{weq3} is locally well-posed in the
Gaussian weighted space $X$. Given $w_0 \in X_0$ with $\|w_0\|_X \le
\epsilon$, let $w \in C^0([0,T_*),X)$ be the maximal solution of
\eqref{weq3} with initial data $w_0$. We know that $w(\tau) \in X_0$
for all $\tau \in [0,T_*)$, and that $w(\tau) \in X_1$ if $w_0 \in
X_1$.  To control the time evolution of $w(\tau)$, we use the energy
estimate
\begin{align*}
  \frac12 \frac{\dd}{\dd \tau} \|w(\tau)\|_X^2  \,&=\,  \langle 
  w(\tau)\,,\,\cL w(\tau)\rangle - \int_{\R^2} G^{-1}(\xi) w(\xi,\tau) 
  v(\xi,\tau) \cdot \nabla w(\xi,\tau)\d\xi \\ \label{energy}
  \,&=\, \langle w(\tau)\,,\,\cL w(\tau)\rangle  +\frac14 
  \int_{\R^2} G^{-1}(\xi) (\xi \cdot v(\xi,\tau)) w(\xi,\tau)^2 \d\xi~,
\end{align*}
where in the first equality we used the skew-symmetry of the operator 
$\Lambda$, and in the second one we integrated by parts and used the 
identity $\nabla G^{-1} = (\xi/2)G^{-1}$. 

Fix $\tau > 0$ and let $f(\xi) = G^{-1/2}(\xi)w(\xi,\tau)$, for all 
$\xi \in \R^2$. Then $\|f\|_{L^2} = \|w(\tau)\|_X$ and, in view 
of \eqref{Lconj},  
\[
  E(f) \,:=\, \|\nabla f\|_{L^2}^2 + \frac{1}{16}\,\| \xi f\|_{L^2}^2
  -\frac12 \|f\|_{L^2}^2 \,=\, -\langle w(\tau)\,,\,\cL 
  w(\tau)\rangle~.
\]
Since $w(\tau) \in X_0$, we have $E(f) \ge \mu\|f\|_{L^2}^2$ 
where $\mu = 1$ if $w_0 \in X_1$ and $\mu = 1/2$ if $w_0 \in 
X_0 \setminus X_1$. In particular, there exists $C > 0$ such that
\[
  E(f) \,\ge\, C(\|\nabla f\|_{L^2}^2 + \|\xi f\|_{L^2}^2 + 
  \|f\|_{L^2}^2)~.
\]
On the other hand, 
\begin{align}\nonumber
  \Big|\int_{\R^2} G^{-1}(\xi) (\xi\cdot v(\xi,\tau)) w(\xi,\tau)^2 
  \d\xi\Big| \,&\le\, \int_{\R^2} |\xi \cdot v(\xi,\tau)| f(\xi)^2 \d\xi \\ 
  \label{trilin} \,&\le\, \|v(\tau)\|_{L^4} \,\|f\|_{L^4} \,\|\xi f\|_{L^2}~.
\end{align}
By \eqref{BSi} we have $\|v\|_{L^4} \le C\|w\|_{L^{4/3}} \le C \|w\|_X$, 
hence the right-hand side of \eqref{trilin} is bounded by 
$C \|w(\tau)\|_X E(f)$. As a consequence, there exists $C_1 > 0$ 
such that
\begin{equation}\label{energy2} 
  \frac12\frac{\dd}{\dd \tau} \|w(\tau)\|_X^2 \,\le\, -E(f)(1 - C_1 
  \|w(\tau)\|_X) \,\le\, -\mu \|w(\tau)\|_X^2 (1 - C_1 
  \|w(\tau)\|_X)~, 
\end{equation}
where the second inequality is valid provided $C_1 \|w(\tau)\|_X \le 1$. 

We now assume that $\epsilon > 0$ is small enough so that $C_1\epsilon 
\le \frac12$. Then the differential inquality \eqref{energy2} implies that 
the norm $\|w(\tau)\|_X$ is nonincreasing with time, hence the solution 
of \eqref{weq3} is global (i.e., $T_* = \infty$) and satisfies 
\eqref{energy2} for all $\tau \ge 0$. Integrating \eqref{energy2}, we obtain 
\[
  \frac{\|w(\tau)\|_X}{1 - C_1 \|w(\tau)\|_X} \,\le\, 
  \frac{\|w_0\|_X}{1 - C_1 \|w_0\|_X}\,e^{-\mu \tau}~, 
  \qquad \hbox{for all }\tau \ge 0~. 
\]
This implies \eqref{wdecay}, since $C_1 \|w_0\|_X \le \frac12$ 
by assumption. 
\QED

Proposition~\ref{locstab} shows that the steady state $w = \alpha G$
of the rescaled vorticity equation \eqref{weq2} is stable with respect
to perturbations in the space $X$, for any value of the circulation
parameter $\alpha \in \R$. Moreover, we have a lower bound on the size
of the local basin of attraction which is uniform in $\alpha$. Here we
mean by ``local basin of attraction'' a neighborhood of the
equilibrium $\alpha G$ which is small enough so that the time
evolution of the perturbations in that neighborhood can be controlled
by the dissipative properties of the linearized operator, as in the
above proof. Of course, Theorem~\ref{thm3} shows that all trajectories
of \eqref{weq2} in $L^1(\R^2)$ converge to the line of equilibria
$\{\alpha G\}_{\alpha \in \R}$ as $\tau \to \infty$, hence one can
argue that the basin of attraction of these equilibria has, actually,
infinite size.

If the initial data are sufficiently localized, Proposition~\ref{locstab} 
allows to improve the conclusion of Theorem~\ref{thm3} by specifying 
sharp convergence rates. Assume for instance that the initial vorticity 
$\omega_0 : \R^2 \to \R$ satisfies
\begin{equation}\label{Gaussbd}
  \int_{\R^2} \omega_0(x)^2 \,\exp\Bigl(\frac{|x|^2}{4\nu t_0}\Bigr)
  \d x \,<\, \infty~,
\end{equation}
for some $t_0 > 0$, and let $\gamma = \int_{\R^2} \omega_0(x)\d x$.  
If $\gamma \neq 0$, we define the center of vorticity as
\[
  x_0 \,=\, \frac{1}{\gamma} \int_{\R^2} x\,\omega_0(x)\d x \,\in\,
  \R^2~,
\]
and we set $x_0 = 0$ otherwise. Using the change of variables
\eqref{wvdef} with $x_0$ as above and $T = t_0$, we transform
\eqref{omeq} into the rescaled vorticity equation \eqref{weq2}, and
hypothesis \eqref{Gaussbd} means that the initial data $w_0$ for
\eqref{weq2} belong to the Gaussian space $X$. Then estimate
\eqref{wupp} (with $\beta > 1/2$) implies that the solution $w(\tau)$
of \eqref{weq2} is uniformly bounded in $X$ for all $\tau \ge 0$, and
arguing as in Section~\ref{ss3.3} it is not difficult to show that
$\|w(\tau) - \alpha G\|_X \to 0$ as $\tau \to \infty$, where $\alpha =
\gamma/\nu$. Thus we can apply Proposition~\ref{locstab} when $\tau$
is sufficiently large, and we conclude that $\|w(\tau) - \alpha G\|_X
= \cO(e^{-\mu \tau})$ as $\tau \to \infty$, for some $\mu \ge 1/2$.
More precisely, when $\alpha \neq 0$, the choice of $x_0$ above
implies that $w_0 - \alpha G \in X_1$, hence we can take $\mu = 1$;
otherwise, we take $\mu = 1/2$. Returning to the original variables,
we obtain in particular the estimate
\begin{equation}\label{Osconv3}
   \Big\|\omega(x,t) - \frac{\gamma}{\nu (t+t_0)} 
   G\Bigl(\frac{x-x_0}{\sqrt{\nu(t+t_0)}}\Bigr)\Big\|_{L^1} \,=\, 
  \cO(t^{-\mu})~, 
  \quad \hbox{as} \quad  t \to \infty~,
\end{equation}
which improves \eqref{Osconv2}. The decay rate in \eqref{Osconv3} 
is optimal in general, but can sometimes be improved by an 
appropriate choice of $t_0$, see \cite{GW2,GR}. We also mention
that the convergence \eqref{Osconv3} still holds if, instead 
of \eqref{Gaussbd}, we assume that the initial vorticity satisfies
$\int_{\R^2} \omega_0(x)^2(1+|x|^2)^m\d x < \infty$ for some $m > 3$, 
see \cite{GW2}. 

To conclude our discussion of Proposition~\ref{locstab}, we would like
to emphasize that the stability analysis of Oseen vortices lead to
conclusions which differ radically from what is known for other
classical examples in fluid mechanics, such as the Poiseuille flow in
a cylindrical pipe or the Couette-Taylor flow between two rotating
cylinders \cite{DR,TE}. Indeed, in the case of Oseen vortices, no
hydrodynamic instability, neither of spectral nor of pseudospectral
nature, develops as the circulation parameter $\alpha = \gamma/\nu$
(which plays the role of the Reynolds number) is increased. On the
contrary, it is possible to show that a fast rotation has even a {\em
  stabilizing effect} on the vortex, as we shall see below.

\subsection{Improved Stability Estimates for Rapidly Rotating Vortices}
\label{ss4.2}

For a more detailed study of the linearized operator $\cL - 
\alpha\Lambda$, it is useful to observe that both operators 
$\cL$ and $\Lambda$ are invariant under rotations about the
origin in $\R^2$. This symmetry is fully exploited if we 
introduce polar coordinates $(r,\theta)$ in $\R^2$ and expand 
the vorticity distribution in Fourier series with respect 
to the angular variable $\theta$. Let us decompose 
\begin{equation}\label{tildeXdef}
  X = \mathop{\oplus}\limits_{n \in \N} \tilde X_n~,
\end{equation}
where $\tilde X_n$ denotes the subspace of all $w \in X$ such 
that 
\[
  w(r\cos\theta,r\sin\theta) = a(r)\cos(n\theta) + b(r)
  \sin(n\theta)~, \qquad r > 0~,\quad \theta \in [0,2\pi]~,
\]
for some radial functions $a, b : \R_+ \to \R$. The invariance 
under rotations implies the following result, which can be 
obtained by a direct calculation\:

\begin{lemma}\label{rotinv} {\rm \cite{GW2}}\\
For each $n \in \N$ the subspace $\tilde X_n$ is invariant under
the action of both linear operators $\cL$ and $\Lambda$. The 
restriction $\cL_n$ of $\cL$ to $\tilde X_n$ is the one-dimensional
operator
\[
   \cL_n \,=\, \partial_r^2 + \Bigl(\frac{r}{2}+\frac{1}{r}
   \Bigr)\partial_r + \Bigl(1 - \frac{n^2}{r^2}\Bigr)~.
\]
The restriction $\Lambda_n$ of $\Lambda$ to $\tilde X_n$ satisfies 
$\Lambda_n \omega \,=\, in(\phi \omega - g\Omega_n)$, where
\begin{equation}\label{phigdef}
   \phi(r) \,=\, \frac{1}{2\pi r^2}(1-e^{-r^2/4})~, \qquad
   g(r) \,=\, \frac{1}{4\pi}\,e^{-r^2/4}~, \qquad r > 0~,
\end{equation}
and (for $n \ge 1$) $\Omega_n$ is the regular solution of 
$-\Omega_n'' -\frac1r \Omega_n' + \frac{n^2}{r^2}\Omega_n = \omega$, 
namely
\begin{equation}\label{Omndef}
   \Omega_n(r) \,=\, \frac{1}{4 n }\left(\int_0^r \Bigl(\frac{r'}{r}
   \Bigr)^n r'\omega(r')\d r' + \int_r^\infty \Bigl(\frac{r}{r'}
   \Bigr)^n r'\omega(r')\d r'\right)~, \quad r > 0~.
\end{equation}
\end{lemma}

In particular, Lemma~\ref{rotinv} implies that the subspace 
$\tilde X_0 \subset X$, which consists of all radially 
symmetric functions, is entirely contained in the kernel 
of the operator $\Lambda$. As was already observed, the 
first order derivatives $\partial_1 G, \partial_2 G$ also 
belong to $\ker(\Lambda)$, and the following result shows
that these are the only non-radially symmetric functions 
in the kernel of $\Lambda$\:

\begin{lemma}\label{kerLam} {\rm \cite{Ma}} $\ker(\Lambda) = \tilde X_0 
\oplus \{\beta_1\partial_1 G + \beta_2 \partial_2 G\,|\, \beta_1,
\beta_2 \in \R\}$.
\end{lemma}

\noindent{\bf Proof.}
We already know that $\tilde X_0 \subset \ker(\Lambda)$. Suppose now
that $w \in \ker(\Lambda) \cap \tilde X_n$ for some $n \ge 1$. 
Without loss of generality, we can assume that $w(r\cos\theta,
r\sin\theta) = \omega(r)\cos(n\theta)$ or $\omega(r)\sin(n\theta)$ for 
some $\omega : \R_+ \to \R$, and using the notations \eqref{phigdef}, 
\eqref{Omndef} we have $\phi \omega - g\Omega_n = 0$. This means that 
$\Omega_n$ satisfies the differential equation
\begin{equation}\label{Omdiff}
  -\Omega_n''(r) -\frac1r \Omega_n'(r) + \Bigl(\frac{n^2}{r^2}
  - \frac{g(r)}{\phi(r)}\Bigr)\Omega_n(r) \,=\, 0~, \qquad 
  r > 0~. 
\end{equation}
Now it is easy to verify that $r^2 g(r)/\phi(r) < 4$ for all $r >
0$. Thus Eq.~\eqref{Omdiff} satisfies the maximum principle if $n \ge
2$, hence has no nontrivial solution vanishing at the origin and at
infinity. This shows that $\ker(\Lambda) \cap \tilde X_n = \{0\}$ for
$n \ge 2$. When $n = 1$, Eq.~\eqref{Omdiff} has a (unique) nontrivial
solution given by $\Omega_1(r) = r\phi(r)$, which satisfies
\eqref{Omndef} with $\omega(r) = rg(r)$. This shows that
$\ker(\Lambda) \cap \tilde X_1$ is spanned by $\partial_1 G$ and
$\partial_2 G$.
\QED

In the rest of this section, to investigate the effect of the
rotation on the stability of Oseen's vortex, we restrict ourselves
to the orthogonal complement of $\ker(\Lambda)$ in $X$. 
This subspace, denoted by $\ker(\Lambda)^\perp$, is invariant 
under the evolution defined by the linearized equation 
$\partial_\tau w = (\cL - \alpha\Lambda)w$, but not under 
the full equation \eqref{weq3}. Let $\cL_\perp$ and $\Lambda_\perp$
denote the restrictions of the operators $\cL$ and $\Lambda$ 
to the invariant subspace $\ker(\Lambda)^\perp$. For any $\alpha 
\in \R$ we define the {\em spectral lower bound}
\[ 
  \Sigma(\alpha) \,=\, \inf\Bigl\{\Re(z)\,\Big|\, z \in 
  \sigma(-\cL_\perp + \alpha\Lambda_\perp)\Bigr\}~,
\]
and the {\em pseudospectral bound} 
\[
  \Psi(\alpha) \,=\, \Bigl(\sup_{\lambda \in \R}\|(\cL_\perp
  -\alpha\Lambda_\perp-i\lambda)^{-1}\|_{X \to X}\Bigr)^{-1}~.
\]
Note that, in the definition of $\Sigma(\alpha)$, we have changed
the sign of the linearized operator, in order to get a positive 
quantity. 

\begin{lemma}\label{SigPsi}
For any $\alpha \in \R$ one has $\Sigma(\alpha) \ge \Psi(\alpha)
\ge 1$. 
\end{lemma}

\noindent{\bf Proof.}
Fix $\alpha \in \R$. Since $\ker(\Lambda)^\perp \subset X_1$, it follows 
from \eqref{sigma1} that $\Sigma(\alpha) \ge 1$. Moreover, if 
$(\cL - \alpha\Lambda + z)w = 0$ for some $z \in \C$ and some 
normalized $w \in \ker(\Lambda)^\perp$, then $(\cL - \alpha\Lambda + 
i\Im(z))w = -\Re(z)w$, so that 
\[
  \Re(z) \,\ge\, \|(\cL_\perp - \alpha\Lambda_\perp +i\Im(z))^{-1}\|^{-1} 
  \,\ge\, \Psi(\alpha)~.
\]
This proves that $\Sigma(\alpha) \ge \Psi(\alpha)$. On the other, 
the proof of Proposition~\ref{specstab} shows that the operator
$\cL_\perp - \alpha\Lambda_\perp +1$ is $m$-dissipative. This in 
particular implies that $\|(\cL_\perp - \alpha\Lambda_\perp-i\lambda
)^{-1}\| \le 1$ for all $\lambda \in \R$, hence $\Psi(\alpha) \ge 1$. 
\QED

Numerical observations due to Prochazka and Pullin \cite{PP} indicate
that $\Sigma(\alpha) = \cO(|\alpha|^{1/2})$ as $|\alpha| \to
\infty$. This has not been proved yet, but a recent result by Maekawa
\cite{Ma} shows that $\Sigma(\alpha) \to +\infty$ as $|\alpha| \to
\infty$. Although the proof is not constructive and does not give any
information on the asymptotic behavior of $\Sigma(\alpha)$, this is
the first rigorous result which shows that a fast rotation stabilizes
Oseen's vortex, in the sense that the decay rate of the perturbations
in $\ker(\Lambda)^\perp$ becomes arbitrarily large as $|\alpha| \to
\infty$. In fact, the argument of \cite{Ma} can be easily modified
to yield the stronger result $\Psi(\alpha) \to \infty$ as 
$|\alpha| \to \infty$. This indicates that the size of the local
basin of attraction of Oseen's vortex $\alpha G$ becomes 
arbitrarily large (in the Gaussian space $X$) when $|\alpha| \to 
\infty$. 

Explicit bounds on $\Sigma(\alpha)$ or $\Psi(\alpha)$ are more
difficult to obtain, and were first established in \cite{GGN} for a
model problem in one space dimension, and in \cite{De1} for a
simplified version of the linearized operator $\cL - \alpha \Lambda$
where the nonlocal term $\Lambda_2$ in $\Lambda$ is omitted.  In both
cases one gets $\Psi(\alpha) = \cO(|\alpha|^{1/3})$ as $|\alpha| \to
\infty$, and for the one-dimensional model we also have $\Sigma(\alpha) 
= \cO(|\alpha|^{1/2})$. These partial results suggest that
$\Sigma(\alpha) \gg \Psi(\alpha)$ when $|\alpha|$ is large, 
a property which reflects the fact that the linearized operator is 
highly non-selfadjoint in the large circulation regime. Indeed, 
it is easily verified that the spectral and pseudospectral bounds 
always coincide for selfadjoint operators.

Although the analysis of the linearized operator $\cL - \alpha
\Lambda$ is still the object of ongoing research, we formulate
the following conjecture which should be proved soon. 

\begin{conjecture}\label{theconj}
There exist constants $\delta > \gamma > 0$ such that
$\Psi(\alpha) = \cO(|\alpha|^\gamma)$ and $\Sigma(\alpha) = 
\cO(|\alpha|^\delta)$ as $|\alpha| \to \infty$. 
\end{conjecture}

On the basis of the analysis of simpler models \cite{GGN}, 
we further conjecture that $\gamma = 1/3$ and $\delta = 1/2$, but 
these optimal values may be difficult to reach. To conclude 
this discussion, we mention that, if the linearized operator $\cL 
- \alpha\Lambda$ is restricted to the invariant subspace
\[
  \mathop{\oplus}\limits_{n \ge N} \tilde X_n \,\subset\, \ker(\Lambda)^\perp~, 
\]
for some sufficiently large integer $N \ge 2$, then it is possible 
to show that the pseudospectral bound in that subspace satisfies 
$\Psi_N(\alpha) =  \cO(|\alpha|^{1/3})$ as $|\alpha| \to \infty$, 
see \cite{De2}. 

\section{Interaction of Vortices in Weakly Viscous Flows}
\label{s5}

Theorem~\ref{thm2} shows that the two-dimensional vorticity equation
\eqref{omeq} is globally well-posed in the space of finite measures,
for any value of the viscosity parameter $\nu$. It is therefore
natural to investigate how the solutions of \eqref{omeq} behave in 
the vanishing viscosity limit. This question is very difficult, because the 
formal limit of \eqref{omeq} as $\nu \to 0$ is the two-dimensional
inviscid vorticity equation
\begin{equation}\label{omeq2}
  \partial_t \omega(x,t) +  u(x,t)\cdot\nabla \omega(x,t) \,=\, 0~,
\end{equation}
which is certainly not well-posed in the space $\cM(\R^2)$. In fact, a
celebrated result of Yudovich \cite{Yu1} shows that \eqref{omeq2} has
a unique global solution if the initial vorticity $\mu$ belongs to
$L^1(\R^2) \cap L^\infty(\R^2)$, see also \cite{Yu2,Vi}. If $\mu \in
\cM(\R^2) \cap H^{-1}(\R^2)$ and if the singular part of $\mu$ has a
definite sign, then \eqref{omeq2} has at least one global weak solution
\cite{De,Maj}, but this result does not cover the case where the 
initial data contain point vortices.

A careful treatment of the inviscid limit is beyond the scope of these
notes, but the following brief discussion will be useful to put our
results into perspective. As a general rule, it is known that the
Navier-Stokes equation is nicely approximated by the Euler equation in
the vanishing viscosity limit if we restrict ourselves to smooth
solutions in a domain without boundary \cite{Sw,Ka3,BM}. Thus the main
difficulties come either from the presence of boundaries or from
singularities in the initial data. Since we chose to work in the whole
plane $\R^2$, we do not mention here the (hard) problems related to
boundary layers, and refer to \cite{SC1,SC2,Gr} for a discussion of
that point. But even in the absence of boundaries, the inviscid limit
may be nontrivial to understand if one considers non-smooth initial
data. In the case of {\em vortex patches}, the situation is well
understood\: the Euler flow is well defined due to Yudovich's theorem,
and approximates the Navier-Stokes flow in every reasonable sense as
$\nu \to 0$, see \cite{CW1,CW2,Ch,Da1,Da2,AD,Hm1,Hm2,Mas,Su}.  The
situation is much less clear for {\em vortex sheets}\: although an
inviscid solution can still be constructed, the Euler flow is unstable
due to the discontinuity of the velocity field (Kelvin-Helmholtz
instability), and the vanishing viscosity limit is as difficult to
study as in the case of Prandtl boundary layers \cite{CS}.  Finally,
the case of {\em point vortices} is somewhat paradoxical\: one the one
hand, these are the most singular initial data for which we can solve
the 2D Navier-Stokes equation, and strictly speaking the Euler flow is
not even defined for such initial data. On the other hand, we have a
nice substitute for the Euler flow in that particular case, namely the
Helmholtz-Kirchhoff point vortex system, which does not exhibit any
dynamical instability. As a result, the inviscid limit can be 
rigorously studied in the case of point vortices. Early results 
in this direction were obtained by Marchioro \cite{Mar1,Mar2}, 
and will be described more precisely below. 

\subsection{The Viscous $N$-Vortex Solution}
\label{ss5.1}

From now on, we fix an integer $N \ge 1$, and we assume that the
initial vorticity is a collection of $N$ point vortices 
characterized by their (pairwise distinct) positions $x_1,\dots,\,x_N 
\in \R^2$ and their (nonzero) circulations $\gamma_1,\dots,\,\gamma_N 
\in \R$\:
\begin{equation}\label{mudef}
  \mu \,=\, \sum_{i=1}^N \gamma_i \,\delta_{x_i}~.
\end{equation}
Given a viscosity $\nu > 0$, we denote by $\omega^\nu(x,t)$,
$u^\nu(x,t)$ the unique solution of the vorticity equation
\eqref{omeq} with initial data $\mu$. Existence of such a solution was
first proved in \cite{BEP}, and uniqueness is guaranteed by
Theorem~\ref{thm2}. Note that Theorem~\ref{thm1} applies if $\sum
|\gamma_i| \le C_0 \nu$, but such an assumption is totally
inappropriate here because we want to study the limit $\nu \to 0$
for fixed initial data. 

If $N = 1$, the solution of \eqref{omeq} is just the Oseen 
vortex with circulation $\gamma_1$ centered at the point $x_1$\:
\[
  \omega^\nu(x,t) \,=\, \frac{\gamma_1}{\nu t}\,G\Bigl(
  \frac{x-x_1}{\sqrt{\nu t}}\Bigr)~,\qquad u^\nu(x,t) \,=\, 
  \frac{\gamma_1}{\sqrt{\nu t}}\,v^G\Bigl(\frac{x-x_1}{\sqrt{\nu 
  t}}\Bigr)~.
\]
When $N \ge 2$, the viscous $N$-vortex solution is no longer explicit, 
but as long as the diffusive length $\sqrt{\nu t}$ is small compared
to the distance between the vortex centers we expect that 
$\omega^\nu(x,t)$ can be approximated by a superposition of 
Oseen vortices\:
\begin{equation}\label{Nvorapp}
  \omega^\nu(x,t) \,\approx\, \sum_{i=1}^N\frac{\gamma_i}{\nu t}\,G
  \Bigl(\frac{x-z_i(t)}{\sqrt{\nu t}}\Bigr)~, \qquad 
  u^\nu(x,t) \,\approx\,  \sum_{i=1}^N \frac{\gamma_i}{\sqrt{\nu t}}\,v^G
  \Bigl(\frac{x-z_i(t)}{\sqrt{\nu t}}\Bigr)~,
\end{equation}
where $z_1(t),\dots,z_N(t)$ denote the positions of the vortex 
centers at time $t$. In the same regime, we also expect that 
the positions $z_i(t)$ will be determined by the Helmholtz-Kirchhoff 
system \cite{He,Ki}\:
\begin{equation}\label{PW}
  z_i'(t) \,=\, \frac{1}{2\pi} \sum_{j\neq i} \gamma_j 
  \,\frac{(z_i(t) - z_j(t))^\perp}{|z_i(t)-z_j(t)|^2}~, 
  \qquad z_i(0) \,=\, x_i~.
\end{equation}
Indeed, we know from the work of Marchioro and Pulvirenti that the 
system of ordinary differential equations \eqref{PW} approximates 
rigorously the motion of localized vortex patches in two-dimensional 
inviscid \cite{Mar0,MP1,MP2} or slightly viscous \cite{Mar1,Mar2} fluids. 

The point vortex system \eqref{PW} is globally well-posed for all
initial data if $N = 2$ or if all circulations $\gamma_i$ have the
same sign. In the general case, however, vortex collisions can occur
in finite time for some exceptional initial configurations
\cite{MP2,Ne}. To eliminate this potential problem, we assume in
what follows that the solution $\{z_1(t),\dots,z_N(t)\}$ of \eqref{PW} is 
defined on some time interval $[0,T]$, and we denote by $d$ 
the minimal distance between the vortex centers\:
\begin{equation}\label{dmin}
  d \,=\, \min_{t \in [0,T]}\,\min_{i\neq j}\,|z_i(t) - z_j(t)| \,>\, 0~.
\end{equation}
We also introduce the turnover time $T_0 = d^2/|\gamma|$, where
$|\gamma| = |\gamma_1| + \dots + |\gamma_N|$, which is a natural time 
scale for the inviscid dynamics described by \eqref{PW}. For instance, 
for a pair of vortices with the same circulation $\gamma$ separated 
by a distance $d$, one can check that the rotation period of each 
vortex around the midpoint is $4\pi^2 T_0$.

With these notations, we can formulate now the main result of 
this section\:

\begin{theorem}\label{thm4} {\rm \cite{Ga}}\\
Assume that the point vortex system \eqref{PW} is well-posed on 
the time interval $[0,T]$. Then the unique solution of the 
two-dimensional vorticity equation \eqref{omeq} with initial data 
$\mu = \sum_{i=1}^N \gamma_i \,\delta_{x_i}$ satisfies
\begin{equation}\label{omnuapp}
  \frac{1}{|\gamma|}\int_{\R^2}\Bigl|\omega^\nu(x,t) - \sum_{i=1}^N 
  \frac{\gamma_i}{\nu t}\,G\Bigl(\frac{x-z_i(t)}{\sqrt{\nu t}}\Bigr)
  \Bigr|\d x \,\le \, K\,\frac{\nu t}{d^2}~, \qquad t \in (0,T]~,
\end{equation}
where $\{z_1(t),\dots,z_N(t)\}$ is the solution of \eqref{PW} and $K$ 
is a dimensionless constant depending only on the ratio $T/T_0$.
\end{theorem}

Theorem~\ref{thm4} describes the viscous $N$-vortex solution to
leading order in our expansion parameter $\sqrt{\nu t}/d$, which is
the ratio of the typical size of the vortex cores to the minimal
distance between the vortex centers. In other words, the approximation
\eqref{Nvorapp} is accurate whenever the vortices are widely separated
compared to their size. Another dimensionless quantity that is 
present in our analysis is the ratio of the observation time $T$ 
to the turnover time $T_0$. It is important to remark that no 
smallness assumption is made on this ratio in Theorem~\ref{thm4}, 
although the constant $K$ in \eqref{omnuapp} becomes very large 
if $T \gg T_0$.  

Since each Oseen vortex converges weakly to a Dirac mass in 
the inviscid limit, an immediate consequence of Theorem~\ref{thm4} 
is\:

\begin{corollary}\label{weakcor}
Under the assumptions of Theorem~\ref{thm4}, the viscous 
$N$-vortex solution $\omega^\nu(x,t)$ satisfies
\begin{equation}\label{weakconv}
  \omega^\nu(\cdot,t) 
  ~\xrightharpoonup[\nu \to 0]{\hbox to 8mm{}}~ 
  \sum_{i=1}^N \gamma_i \,\delta_{z_i(t)}~, \quad 
  \hbox{for all } t \in [0,T]~,
\end{equation}
where $\{z_1(t),\dots,\,z_N(t)\}$ is the solution of \eqref{PW}. 
\end{corollary}

While weaker than Theorem~\ref{thm4}, this corollary provides a
rigorous derivation of the point vortex system \eqref{PW} using the
Navier-Stokes equation in the inviscid limit. In contrast, the
classical approach of Marchioro and Pulvirenti \cite{MP1,MP2} 
allows to justify \eqref{PW} within the framework of Euler's 
equation, but requires an approximation argument which is
not needed in Corollary~\ref{weakcor}.

\subsection{Decomposition into Vorticity Profiles}
\label{ss5.2}

In the rest of this section, we give a sketch of the 
proof of Theorem~\ref{thm4}. Our starting point is a natural
decomposition of the viscous $N$-vortex solution, which follows
immediately from the representation formula \eqref{representation}. 
Since the initial measure $\mu$ is a sum of Dirac masses, we
have by \eqref{representation}
\begin{equation}\label{omdecomp}
  \omega^\nu(x,t) \,=\, \sum_{i=1}^N \omega_i^\nu(x,t)~, \qquad
  u^\nu(x,t) \,=\, \sum_{i=1}^N u_i^\nu(x,t)~, 
\end{equation}
where
\[
  \omega_i^\nu(x,t) \,=\, \gamma_i \Gamma_u^\nu(x,t;x_i,0)~, 
  \qquad x \in \R^2~, \quad 0 < t \le T~, \quad i = 1,\dots,N~,
\]
and $u_i^\nu(x,t)$ is the velocity field corresponding to
$\omega^\nu(x,t)$ via the Biot-Savart law \eqref{BS}. Here
$\Gamma_u^\nu(x,t;y,s)$ is the fundamental solution of the
convection-diffusion equation \eqref{Ueq} with advection field $U(x,t)
= u^\nu(x,t)$. In view of \eqref{gamm1}, we have the bound
\begin{equation}\label{omegaibdd}
   |\omega_i^\nu(x,t)| \,\le\, K_1\,\frac{|\gamma_i|}{\nu t} 
   \,\exp\Bigl(-\beta\frac{|x-x_i|^2}{4\nu t}\Bigr)~,   
\end{equation}
for all $i \in \{1,\dots,N\}$, all $x \in \R^2$, and all $t \in
(0,T]$, where $\beta \in (0,1)$ and the dimensionless constant 
$K_1 > 0$ depends only on the total circulation $|\gamma|$ and on the
viscosity $\nu$. Estimate \eqref{omegaibdd} gives a precise
information on the $N$-vortex solution when $\nu > 0$ is fixed and $t
\to 0$, but cannot be used to control $\omega_i^\nu(x,t)$ when $t > 0$
is fixed and $\nu \to 0$, because the constant $K_1$ blows up in that
regime.

To analyze more precisely each term in the decomposition 
\eqref{omdecomp} we introduce, for each $i \in \{1,\dots,N\}$, 
the self-similar variable
\[
  \xi \,=\, \frac{x - z_i(t)}{\sqrt{\nu t}}~, 
\]
where $\{z_1(t),\dots,\,z_N(t)\}$ is the solution of \eqref{PW}.
As in \eqref{wvdef}, we define rescaled vorticities $w_i^\nu(\xi,t) 
\in \R$ and velocities $v_i^\nu(\xi,t) \in \R^2$ by setting
\begin{equation}\label{wvdef2}
  \omega_i^\nu(x,t) \,=\, \frac{\gamma_i}{\nu t}\,w_i^\nu
  \Bigl(\frac{x - z_i(t)}{\sqrt{\nu t}}\,,\,t\Bigr)~, \qquad
  u_i^\nu(x,t) \,=\, \frac{\gamma_i}{\sqrt{\nu t}}\,v_i^\nu
  \Bigl(\frac{x - z_i(t)}{\sqrt{\nu t}}\,,\,t\Bigr)~.
\end{equation}
Our goal is to describe the vorticity profiles $w_i^\nu(\xi,t)$ as
precisely as possible, in the regime where $\sqrt{\nu t} \ll d$.
According to \eqref{omnuapp}, we have $w_i^\nu(\xi,t) = G(\xi) +
\cO(\nu t/d^2)$ for each $i \in \{1,\dots,N\}$, where $G$ is the
profile of Oseen's vortex. However, it is very important to realize
that this leading order approximation is {\em not sufficient} to
control the inviscid limit of the viscous $N$-vortex solution. The
reason is that the profile $G$ is radially symmetric, and thus does
not take into account the deformations of the vortices due to mutual
interactions.

If we go back to the explicit formula \eqref{Osdef}, we see that 
the velocity field of Oseen's vortex is very large near the center if
the viscosity $\nu$ is small, with a maximal angular speed of the
order of $|\gamma|/(\nu t)$. As long as the vortex stays isolated,
this large velocity does not affect the vorticity distribution, which 
is perfectly radially symmetric. Now, if the same vortex is advected 
by a non-homogenous external field, which in our case will be the
velocity field produced by the other $N{-}1$ vortices, its vorticity 
profile gets deformed under the external strain and the vortex 
starts feeling the effect of its own velocity field. As is easily verified, 
this {\em self-interaction} effect is very strong if the circulation
Reynolds number $|\gamma|/\nu$ is large, even for a moderate
deformation of the vortex core. This gives a clear indication that 
self-interaction effects must be taken into account if one wants to 
control the vanishing viscosity limit in presence of point vortices. 

The deformation of a rapidly rotating Oseen vortex advected by an 
external velocity field can be computed using the properties of 
the linearized operator $\cL - \alpha\Lambda$, which was defined
in Section~\ref{s4}. Remarkably enough, one finds that Oseen's 
vortex adapts its shape in such a way that the self-interaction 
{\em counterbalances} the effect of the external field \cite{TT,TK}, 
except for radially symmetric corrections and for a rigid translation. 
This is related to Lemma~\ref{kerLam}, which shows that the kernel
of the operator $\Lambda$ consists of all radially symmetric functions
and of the two-dimensional subspace spanned by $\partial_1 G, 
\partial_2 G$. This fundamental fact explains why one can
observe, in turbulent two-dimensional flows, stable asymmetric
vortices which in a first approximation are simply advected by the
main stream \cite{BMMPW}. 

In Section~\ref{ss5.3}, we implement this idea by constructing 
a higher order approximation of the viscous $N$-vortex solution 
which takes into account the self-interaction effects. For each 
$i \in \{1,\dots,N\}$ and all $t \in [0,T]$, the approximated vortex 
profile will be of the form
\begin{equation}\label{wapp}
  w_i^\app(\xi,t) \,=\, G(\xi) + \Bigl(\frac{\nu t}{d^2}\Bigr)
  \Bigl\{\bar F_i(\xi,t) + F_i^\nu(\xi,t)\Bigr\}~, \qquad 
  \xi \in \R^2~,
\end{equation}
where $\bar F_i(\xi,t)$ is a radially symmetric function of 
$\xi$ which we will not describe precisely because it 
only represents a small perturbation of the Gaussian profile
$G(\xi)$. The important correction at this order is the 
nonsymmetric function $F_i^\nu(\xi,t)$ which is given by
\begin{equation}\label{Finuexp}
  F_i^\nu(\xi,t) \,=\, \frac{d^2}{4\pi}\,\omega(|\xi|) 
  \sum_{j\neq i} \frac{\gamma_j}{\gamma_i}\,\frac{1}{|z_{ij}(t)|^2}
  \Bigl(2\frac{|\xi \cdot z_{ij}(t)|^2}{|\xi|^2 |z_{ij}(t)|^2}
  - 1\Bigr) + \cO\Bigl(\frac{\nu}{|\gamma|}\Bigr)~, 
\end{equation}
where we denote $z_{ij}(t) = z_i(t) - z_j(t)$. Here $\omega : (0,\infty)
\to \R$ is a smooth, positive function satisfying $\omega(r)
\approx C_1 r^2$ as $r \to 0$ and $\omega(r) \approx C_2 r^4
e^{-r^2/4}$ as $r \to \infty$ for some $C_1, C_2 > 0$. 

The right-hand side of \eqref{wapp} is the beginning of an asymptotic
expansion of the rescaled vortex patch $w_i^\nu(\xi,t)$ in powers of
the non-dimensional parameter $\sqrt{\nu t}/d$, which is the ratio of
the size of the vortex cores to the minimal distance between the
centers. Each term in this expansion can in turn be developed in
powers of the inverse circulation Reynolds number $\nu/|\gamma|$.  The
most important effect is due to the nonsymmetric term
$F_i^\nu(\xi,t)$, which describes to leading order the deformation of
the $i^{\rm th}$ vortex due to the influence of the other
vortices. Keeping only that term and using polar coordinates $\xi =
(r\cos\theta,r\sin\theta)$, we can rewrite \eqref{wapp} in the
following simplified form
\begin{equation}\label{wapp2}
  w_i^\app(\xi,t) \,=\, g(r) + \frac{\omega(r)}{4\pi} \sum_{j\neq i} 
  \frac{\gamma_j}{\gamma_i}\,\frac{\nu t}{|z_{ij}(t)|^2}
  \cos\Bigl(2(\theta - \theta_{ij}(t))\Bigr) \,+\, \dots~,
\end{equation}
where $g(|\xi|) = G(\xi)$ and $\theta_{ij}(t)$ is the argument of the 
planar vector $z_{ij}(t) = z_i(t) - z_j(t)$. This formula 
allows us to compute the principal axes and the eccentricities of 
the vorticity contours, which are elliptical at this level of 
approximation. 

To formulate our result in its final form, we introduce a function 
space that is appropriate to control the vorticity profiles. 
Given a small $\beta \in (0,1)$, which will be specified later, we define 
the weighted $L^2$ space $X_\beta$ equipped with the norm
\begin{equation}\label{Xdef}
  \|w\|_{X_\beta} \,=\, \Bigl(\int_{\R^2} |w(\xi)|^2 \,e^{\beta |\xi|/4}
  \d\xi\Bigr)^{1/2}~.
\end{equation}
In particular, we have $X_\beta \hookrightarrow L^1(\R^2)$ for 
any $\beta > 0$. As is explained in \cite{Ga}, it would be more 
natural to use here the Gaussian space $X$ defined in 
\eqref{Xscalar}, but for technical reasons it appears necessary 
to replace the Gaussian weight by an exponential one, at least
if the ratio $T/T_0$ is large. 

\begin{proposition}\label{propfinal}
Assume that the point vortex system \eqref{PW} is well-posed on 
the time interval $[0,T]$, and let $\omega^\nu(x,t)$ be the 
solution of \eqref{omeq} with initial data \eqref{mudef}. 
Then there exist positive constants $\beta$ and $K_2$, 
depending only on the ratio $T/T_0$, such that, if $\omega^\nu(x,t)$ 
is decomposed as in \eqref{omdecomp}, then the rescaled profiles 
$w_i^\nu(\xi,t)$ defined by \eqref{wvdef2} satisfy
\begin{equation}\label{strongconv}
  \max_{i=1,\dots,N} \|w_i^\nu(\cdot,t) - w_i^\app(\cdot,t)\|_{X_\beta} 
  \,\le\, K_2 \Bigl(\frac{\nu t}{d^2}\Bigr)^{3/2}~, 
\end{equation}
for all $t \in (0,T]$, provided $\nu T/d^2 \le K_2^{-1}$. 
\end{proposition}

As is clear from \eqref{strongconv}, \eqref{wapp}, and \eqref{wvdef2},
Proposition~\ref{propfinal} implies immediately Theorem~\ref{thm4}.
Note that the error term on the right-hand side of \eqref{strongconv}
is smaller than in \eqref{omnuapp}, and in particular smaller than the
first order corrections to the Gaussian profile in \eqref{wapp}.  This
means that the deformations of the interacting vortices are indeed
described, to leading order, by \eqref{wapp2}. According to that
formula, each vortex adapts its shape {\em instantaneously} to the
relative positions of the other vortices, without oscillations or
inertia.  By this we mean that, for each $t \in (0,T]$, the angular
factor $\cos(2(\theta - \theta_{ij}))$ in \eqref{wapp2}, which gives
the leading order deformation up to a time-dependent prefactor $\nu
t/|z_{ij}|^2$, is entirely determined by the instantaneous positions
of the vortex centers. In contrast, it is shown in \cite{Ga} that the
first order radially symmetric corrections $\bar F_i(\xi,t)$ do not
only depend on the instantaneous vortex positions, but on the whole
history of the system.

It is instructive at this point to recall what is known in the more
general situation where the initial data are not point vortices, but
finite size vortex patches. For simplicity, we consider the case of
two identical, radially symmetric vortex patches of size $R > 0$
initially located at points $x_1, x_2 \in \R^2$ with $|x_1 - x_2| \gg
R$. Under the time evolution defined by \eqref{omeq}, the centers of
the patches will rotate with approximatively constant angular velocity
around the mid-point $(x_1 + x_2)/2$. To obtain a more precise
description of the solution, we go to a rotating frame where the
centers of the patches stay fixed, and we concentrate on the
deformations of the vortex cores. What is observed in numerical
simulations \cite{DV} is that the interaction begins with a fast
relaxation process, during which each vortex adapts its shape to the
velocity field generated by the other vortex. This first step depends
on the details of the initial data, and is characterized by temporal
oscillations of the vortex cores which disappear on a non-viscous time
scale. In a second step, the vortices relax to a Gaussian-like profile
at a diffusive rate, and the system reaches a ``metastable state''
which is independent of the initial data, and will persist until two
vortices get sufficiently close to start a merging process. In this
metastable regime, the vortex centers move in the plane according to
the Helmholtz-Kirchhoff dynamics, and the vortex profiles are uniquely
determined, up to a scaling factor, by the relative positions of the
centers. This is exactly the situation described by
Proposition~\ref{propfinal}. In fact, if we start with point vortices,
the system is immediately in the metastable state which, for more
general initial data, is reached only after the transients steps
described above. In this sense, point vortices can be considered as
{\em well-prepared initial data} for the vortex interaction problem.

\subsection{Perturbation Expansion and Error Estimates}
\label{ss5.3}

In this final section, we briefly indicate how to construct an 
asymptotic expansion of the viscous $N$-vortex solution and
to control the error terms. We first write the evolution equation 
satisfied by the rescaled vorticity profiles defined in 
\eqref{wvdef2}. Since $\partial_t \omega_i^\nu + u^\nu\cdot \nabla
\omega_i^\nu = \nu \Delta \omega_i^\nu$ for $i = 1,\dots,N$, 
where $u^\nu = \sum_{j=1}^N u_j^\nu$, we obtain
\[ 
  t\partial_t w_i^\nu(\xi,t) + \left\{\sum_{j=1}^N \frac{\gamma_j}{\nu}
  \,v_j^\nu\Bigl(\xi + \frac{z_{ij}(t)}{\sqrt{\nu t}}\,,\,t\Bigr)
  - \sqrt{\frac{t}{\nu}}\,z_i'(t)\right\}\cdot \nabla w_i^\nu(\xi,t) 
  \,=\, (\cL w_i^\nu)(\xi,t)~,
\]
where $\cL$ is the linear operator \eqref{weq2}. The left-hand side
is clearly singular in the limit $\nu \to 0$, but this difficulty
can be partially eliminated by an appropriate choice of the speeds 
$z_i'(t)$. To this purpose, we set
\begin{equation}\label{PW2}
  z_i'(t) \,=\, \sum_{j=1}^N \frac{\gamma_j}{\sqrt{\nu t}}
  \, v^G\Bigl(\frac{z_{ij}(t)}{\sqrt{\nu t}}\Bigr)~, \qquad
  i = 1,\dots,N~,
\end{equation}
where $z_{ij}(t) = z_i(t) - z_j(t)$ and $v^G$ is the velocity 
profile defined in \eqref{Gdef}. This is a viscous approximation
of the original point vortex system \eqref{PW}, which has its own
interest and was recently studied in \cite{NSUW,JKN}. It is shown in
\cite[Lemma~2]{Ga} that the solutions of \eqref{PW2} stay extremely
close to those of \eqref{PW} in our perturbative regime where the
vortices are widely separated from each other. In what follows, we
thus make no distinction between the solutions of \eqref{PW} and
\eqref{PW2}. 
 
As was already mentioned, we expect the Gaussian function $G$ 
to represent the leading order approximation of each vorticity profile. 
Substituting $G$ for $w_i^\nu(\xi,t)$ in the evolution equation
above, and using \eqref{PW2} together with $\partial_t G = \cL G 
= 0$, we obtain the following expression for the residuum of this 
naive approximation\:
\[
  R_i^{(0)}(\xi,t) \,=\, \sum_{j\neq i} \frac{\gamma_j}{\nu}
  \biggl\{v^G\Bigl(\xi + \frac{z_{ij}(t)}{\sqrt{\nu t}}\Bigr)
    - v^G\Bigl(\frac{z_{ij}(t)}{\sqrt{\nu t}}\Bigr)\biggr\}\cdot 
  \nabla G(\xi)~. 
\]
This expression is not singular in the limit $\nu \to 0$, and 
a direct calculation (see \cite[Proposition 1]{Ga}) yields the 
asymptotic expansion
\begin{equation}\label{R0exp}
  R_i^{(0)}(\xi,t) \,=\, \frac{\gamma_i t}{d^2}\biggl\{
  A_i(\xi,t) + \Bigl(\frac{\nu t}{d^2}\Bigr)^{1/2}
  B_i(\xi,t) +  \Bigl(\frac{\nu t}{d^2}\Bigr) C_i(\xi,t) 
  + \tilde R_i^{(0)}(\xi,t)\biggr\}~,
\end{equation}
for all $\xi \in \R^2$ and all $t \in (0,T]$, where
\begin{align*}
  A_i(\xi,t) \,&=\, \frac{d^2}{2\pi} \sum_{j \neq i}
  \frac{\gamma_j}{\gamma_i}\,\frac{(\xi\cdot z_{ij}(t))
  (\xi\cdot z_{ij}(t)^\perp)}{|z_{ij}(t)|^4}\,G(\xi)~,\\ 
  B_i(\xi,t) \,&=\, \frac{d^3}{4\pi} \sum_{j \neq i}
  \frac{\gamma_j}{\gamma_i}\,\frac{(\xi\cdot z_{ij}(t)^\perp)}
  {|z_{ij}(t)|^6}\,\Bigl(|\xi|^2 |z_{ij}(t)|^2 - 4 
  (\xi\cdot z_{ij}(t))^2\Bigr)\,G(\xi)~,\\ 
  C_i(\xi,t) \,&=\, \frac{d^4}{\pi} \sum_{j \neq i}
  \frac{\gamma_j}{\gamma_i}\,\frac{(\xi\cdot z_{ij}(t))
  (\xi\cdot z_{ij}(t)^\perp)}{|z_{ij}(t)|^8}\Bigl(
  2(\xi\cdot z_{ij}(t))^2 - |\xi|^2 |z_{ij}(t)|^2\Bigr)\,G(\xi)~.
\end{align*}
Moreover, given any $\delta < 1$, the last term $\tilde R_i^{(0)}$ in 
\eqref{R0exp} can be estimated as follows\:
\[
   |\tilde R_i^{(0)}(\xi,t)| \,\le\, C\Bigl(\frac{\nu t}{d^2}
   \Bigr)^{3/2}\,e^{-\delta|\xi|^2/4}~, \qquad \xi \in \R^2~, 
   \quad 0 < t \le T~. 
\]

The formula \eqref{R0exp} shows in particular that $R_i^{(0)}(\xi,t) =
\cO(1)$ as $\nu \to 0$. Thus, if we decompose $w_i^\nu(\xi,t) = G(\xi)
+ \tilde w_i(\xi,t)$, the equation for $\tilde w_i(\xi,t)$ will
contain a source term of size $\cO(1)$ as $\nu \to 0$, and we
therefore expect that the remainder $\tilde w_i(\xi,t)$ itself will be
of size $\cO(1)$ after a short time.  But, as is easily verified, the
equation for $\tilde w_i(\xi,t)$ contains nonlinear terms involving
negative powers of $\nu$, and such terms cannot be controlled in the
vanishing viscosity limit if $\tilde w_i(\xi,t)$ is $\cO(1)$. This is
the reason why it is necessary to construct a more precise approximate
solution of our vorticity profiles, in order to ``desingularize'' the
equation for the remainder.

To this end, we first consider approximate solutions of the form
\[
  w_i^\app(\xi,t) \,=\, G(\xi) + \Bigl(\frac{\nu t}{d^2}\Bigr)
  F_i(\xi,t)~, \qquad
  v_i^\app(\xi,t) \,=\,  v^G(\xi) + \Bigl(\frac{\nu t}{d^2}\Bigr)
  v^{F_i}(\xi,t)~,
\]
where $F_1(\xi,t),\dots, F_n(\xi,t)$ are smooth vorticity profiles
to be determined, and the corresponding velocity profiles 
$v^{F_1}(\xi,t),\dots,v^{F_n}(\xi,t)$ are obtained via the 
Biot-Savart law \eqref{BS}. The residuum of this improved approximation 
can be computed \cite{Ga}, and has the form
\[
  R_i^{(1)}(\xi,t) \,=\, \frac{\gamma_i t}{d^2}
  \Bigl(A_i(\xi,t) + \Lambda F_i(\xi,t)\Bigr) + \cO\left(
  \frac{\sqrt{\nu t}}{d}\right)~,
\]
where $\Lambda$ is the linear operator defined in \eqref{Lamdef}.  To
minimize the residuum, we want to choose $F_i$ so as to satisfy the
linear equation $\Lambda F_i + A_i = 0$. This is indeed possible,
because $A_i$ belongs to the Schwartz class $\cS(\R^2)$ and to the
subspace $\tilde X_2 \subset X$, see \eqref{tildeXdef}, and it is
shown in \cite{MKO,GW3,Ga} that $\cS(\R^2) \cap \tilde X_2 \subset
\ran(\Lambda)$. The equation $\Lambda F_i + A_i = 0$ has therefore a
unique solution in $\tilde X_2$, which coincides with the leading term
in the expression \eqref{Finuexp} of $F_i^\nu(\xi,t)$.  Note however
that the equation above determines $F_i$ only up to an element of
$\ker(\Lambda)$. As a matter of fact, in \eqref{wapp}, a radially
symmetric correction $\bar F_i(\xi,t)$ is added to $F_i^\nu(\xi,t)$ in
order to compensate for higher order terms (which will not be
considered here). In any case, we can choose $F_i$ in such a way that
$R_i^{(1)}(\xi,t) = \cO(\sqrt{\nu t}/d)$.

Following the same procedure, we next construct a higher order 
approximate solution of the form
\begin{align}\label{wiapp2}
  w_i^\app(\xi,t) \,&=\,  {\DS G(\xi) + \Bigl(\frac{\nu t}{d^2}\Bigr)
  F_i(\xi,t) + \Bigl(\frac{\nu t}{d^2}\Bigr)^{3/2}H_i(\xi,t) + 
  \Bigl(\frac{\nu t}{d^2}\Bigr)^2 K_i(\xi,t)~,} \\ \nonumber
  v_i^\app(\xi,t) \,&=\, {\DS v^G(\xi) + \Bigl(\frac{\nu t}{d^2}\Bigr)
  v^{F_i}(\xi,t) + \Bigl(\frac{\nu t}{d^2}\Bigr)^{3/2}v^{H_i}(\xi,t) + 
  \Bigl(\frac{\nu t}{d^2}\Bigr)^2 v^{K_i}(\xi,t)~.}
\end{align}
Again, the vorticity profiles $F_i$, $H_i$, and $K_i$ are obtained
by solving under-determined linear equations involving the operator
$\Lambda$. More precisely, we have $\Lambda F_i + A_i = 0$ and 
$\Lambda H_i + B_i = 0$, but the equation for $K_i$ is more 
complicated. The residuum of that approximation now satisfies
\begin{equation}\label{Ri3}
  |R_i^{(3)}(\xi,t)| \,\le\,  C \,\Bigl(\frac{\nu t}{d^2}\Bigr)^{3/2}
  \,e^{-\delta|\xi|^2/4}~,  \qquad \xi \in \R^2~, 
   \quad 0 < t \le T~. 
\end{equation}

Finally, once the approximate solution \eqref{wiapp2} has been
constructed, we decompose the vorticity profiles as follows\:
\begin{equation}\label{tildewidef}
  w_i^\nu(\xi,t) \,=\, w_i^\app(\xi,t) + \Bigl(\frac{\nu t}{d^2}\Bigr)
  \,{\tilde w}_i(\xi,t)~, \quad
  v_i^\nu(\xi,t) \,=\, v_i^\app(\xi,t) + \Bigl(\frac{\nu t}{d^2}\Bigr)
  \,{\tilde v}_i(\xi,t)~,
\end{equation}
and we consider the evolution system satisfied by the remainder 
$\tilde w_i(\xi,t)$, $\tilde v_i(\xi,t)$. This system reads
\begin{align}\nonumber
 t\partial_t \tilde w_i(\xi,t) &- (\cL \tilde w_i)(\xi,t) 
  + \tilde w_i(\xi,t) \\[2mm]\nonumber 
 & + \frac{\gamma_i}{\nu} \Bigl(v_i^\app(\xi,t)\cdot \nabla
   \tilde w_i(\xi,t) + \tilde v_i(\xi,t)\cdot \nabla w_i^\app(\xi,t)
   \Bigr) \\[1mm] \nonumber
 & + \sum_{j\neq i} \frac{\gamma_j}{\nu}\left\{v_j^\app
   \Bigl(\xi + \frac{z_{ij}(t)}{\sqrt{\nu t}}\,,\,t\Bigr)
   - v^G\Bigl(\frac{z_{ij}(t)}{\sqrt{\nu t}}\Bigr)\right\}\cdot 
   \nabla \tilde w_i(\xi,t) \\ \label{finsys}
 & + \sum_{j\neq i} \frac{\gamma_j}{\nu}\,\tilde v_j
   \Bigl(\xi + \frac{z_{ij}(t)}{\sqrt{\nu t}}\,,\,t\Bigr)
   \cdot \nabla w_i^\app(\xi,t) \\ \nonumber
 & + \sum_{j=1}^N \frac{\gamma_j t}{d^2}\,\tilde v_j
   \Bigl(\xi + \frac{z_{ij}(t)}{\sqrt{\nu t}}\,,\,t\Bigr)
   \cdot \nabla \tilde w_i(\xi,t) + \tilde R_i(\xi,t) \,=\, 0~,  
\end{align}
and using \eqref{Ri3}, \eqref{tildewidef} we see that the residuum
$\tilde R_i$ is $\cO(\sqrt{\nu t}/d)$ as $\nu \to 0$. Also, due to the
ansatz \eqref{tildewidef}, the nonlinear terms in \eqref{finsys} are
now regular in the limit $\nu \to 0$. The linear terms still contain
negative powers of $\nu$, but this difficulty can be avoided by using
appropriate norms. In a much simpler setting, the same idea was
already present in the proof of Proposition~\ref{locstab}, where the
skew-symmetry of the operator $\Lambda$ in the Gaussian space $X$ 
was used to obtain a stability estimate that was uniform in the 
circulation parameter $\alpha$. 

To control the solution of \eqref{finsys}, we introduce a weighted 
energy functional of the form 
\[
  E(t) \,=\, \sum_{i=1}^N \int_{\R^2} p_i(\xi,t) 
  |\tilde w_i(\xi,t)|^2 \d\xi~, \qquad t \in (0,T]~.
\]
If the observation time $T > 0$ is small with respect to the turnover 
time $T_0 = d^2/|\gamma|$, we can take $p_i(\xi,t) = p_{a(t)}(\xi)$ 
for $i = 1,\dots,N$, where $a(t) = d/(3\sqrt{\nu t})$ and
\[
  p_a(\xi) \,=\, \left\{\begin{array}{ccl}
  e^{|\xi|^2/4} & \hbox{if} & |\xi| \le a~, \\
  e^{a^2/4} & \hbox{if} & a \le |\xi| \le Ka~, \\
  e^{|\xi|^2/(4K^2)} & \hbox{if} & |\xi| \ge Ka~,
  \end{array}\right.
\]
for some $K \gg 1$. We then have $e^{|\xi|^2/(4K^2)} \le p_i(\xi,t) 
\le e^{|\xi|^2/4}$ for all $\xi \in \R^2$ and $t \in (0,T)$. 
With this choice, we obtain from \eqref{finsys} a differential 
inequality for the weighted energy $E(t)$, which can be integrated using 
Gronwall's lemma and yields the bound\:
\[
  \int_{\R^2} e^{\frac{|\xi|^2}{4K^2}} \Bigl(|\tilde w_1(\xi,t)|^2 
  + \dots + |\tilde w_N(\xi,t)|^2 \Bigr)\d\xi \,\le\, E(t) \,\le\,
  C\,\frac{\nu t}{d^2}~.
\]
This concludes the proof of Proposition~\ref{propfinal} if $T \ll T_0$.  
In the general case, one has to introduce more complicated weights, 
which can be constructed using the same procedure as the approximate
solution itself. These weights satisfy $e^{\beta|\xi|/4} \le p_i(\xi,t) 
\le e^{|\xi|^2/4}$, for some small $\beta > 0$ depending only on $T/T_0$. 
We thus obtain the weaker estimate\:
$$
  \int_{\R^2} e^{\frac{\beta |\xi|}{4}} \Bigl(|\tilde w_1(\xi,t)|^2 
  + \dots + |\tilde w_N(\xi,t)|^2 \Bigr)\d\xi \,\le\, E(t) \,\le\,
  C\,\frac{\nu t}{d^2}~,
$$
which still implies the desired conclusion. 

\medskip\noindent{\bf Acknowledgements.} 
I would like to thank Josef Malek and Mirko Rokyta for inviting me
to give a series of lectures at the (perfectly organized) school
entitled "Mathematical Theory in Fluid Mechanics", which was held in
Kacov in May-June 2011. The present notes contain a slightly expanded
version of the content of my course. They also benefited from
previous lectures I had the opportunity to deliver on the same topics
at the 1st Franco-Brazilian Fluids Summer School in Campinas (Brazil)
in January 2010, and at a meeting on fluid mechanics in Etretat
(France) in September 2010.


\begin{thebibliography}{99}
\setlength{\itemsep}{-0.4mm}

\bibitem{AD}
\newblock H. Abidi and R. Danchin,
\newblock \emph{Optimal bounds for the inviscid limit of Navier-Stokes 
  equations}, 
  Asymptot. Anal. \textbf{38} (2004), 35--46. 

\bibitem{BM}
\newblock Th. Beale and A. Majda, 
\newblock \emph{Rates of convergence for viscous splitting of the 
  Navier-Stokes equations},
  Math. Comp. \textbf{37} (1981), 243--259

\bibitem{BA}
\newblock  M. Ben-Artzi,
\newblock \emph{Global solutions of two-dimensional Navier-Stokes and 
Euler equations}, 
Arch. Rational Mech. Anal., \textbf{128} (1994), 329--358. 

\bibitem{BEP}
\newblock G. Benfatto, R. Esposito and M. Pulvirenti, 
\newblock \emph{Planar Navier-Stokes flow for singular initial data},
Nonlinear Anal., \textbf{9} (1985), 533--545. 

\bibitem{BMMPW} 
\newblock A. Bracco, J.C. McWilliams, G. Murante, A. Provenzale, 
and J.B. Weiss,
\newblock \emph{Revisiting freely-decaying two-dimensional turbulence 
at millennial resolution}, 
Physics of Fluids \textbf{12} (2000), 2931--2941. 

\bibitem{Br}
\newblock H. Brezis, 
\newblock \emph{Remarks on the preceding paper by M. Ben-Artzi: 
``Global solutions of two-dimensional Navier-Stokes and Euler 
equations''}, 
Arch. Rational Mech. Anal., \textbf{128} (1994), 359--360. 

\bibitem{CS}
\newblock R. Caflisch and M. Sammartino,
\newblock \emph{Vortex layers in the small viscosity limit},
  ``WASCOM 2005''--13th Conference on Waves and Stability in 
  Continuous Media, 59--70, World Sci. Publ., Hackensack, NJ, 2006.

\bibitem{CPR}
\newblock E. Caglioti, M. Pulvirenti and F. Rousset,
\newblock \emph{On a constrained 2-D Navier-Stokes equation}, 
Commun. Math. Phys. \textbf{290} (2009), 651–-677.

\bibitem{CP}
\newblock M. Cannone, M. and F. Planchon,
\newblock \emph{Self-similar solutions for Navier-Stokes equations in 
$\R^3$}, Commun. Partial Differential Equations \textbf{21} (1996),  
179–-193.

\bibitem{CL}
\newblock E. A. Carlen and M. Loss,
\newblock  \emph{Optimal smoothing and decay estimates for viscously 
 damped conservation laws, with applications to the $2$-D Navier-Stokes
  equation},
Duke Math. J., \textbf{81} (1995), 135--157 (1996). 

\bibitem{Ch}
\newblock J.-Y. Chemin,
\newblock \emph{A remark on the inviscid limit for two-dimensional 
  incompressible fluids},
  Comm. Partial Diff. Equations \textbf{21} (1996), 1771--1779.

\bibitem{CDGG}
\newblock J.-Y. Chemin, B. Desjardins, I. Gallagher and E. Grenier, 
\newblock "Mathematical Geophysics. An Introduction to Rotating Fluids 
and the Navier-Stokes Equations", Oxford Lecture Series in Mathematics 
and its Applications \textbf{32},  Oxford University Press, 2006. 

\bibitem{Con}
\newblock P. Constantin,
\newblock \emph{On the Euler equations of incompressible fluids}, 
Bull. Amer. Math. Soc. \textbf{44} (2007), 603--621.

\bibitem{CW1}
\newblock P. Constantin and J. Wu,
\newblock \emph{Inviscid limit for vortex patches},
  Nonlinearity \textbf{8} (1995), 735--742. 

\bibitem{CW2}
\newblock P. Constantin and J. Wu, 
\newblock \emph{The inviscid limit for non-smooth vorticity}, 
  Indiana Univ. Math. J. \textbf{45} (1996), 67--81.

\bibitem{Co}
\newblock G.-H. Cottet,
\newblock \emph{\'Equations de Navier-Stokes dans le plan avec tourbillon
  initial mesure},
C. R. Acad. Sci. Paris S\'er. I Math., \textbf{303} (1986),  
105--108. 

\bibitem{Cou} 
\newblock Y. Couder,
\newblock \emph{Observation exp\'erimentale de la turbulence 
bidimensionnelle dans un film liquide mince}, 
C. R. Acad. Sci. Paris II \textbf{297} (1983), 641--645. 

\bibitem{Da1}
\newblock R. Danchin,
\newblock \emph{Poches de tourbillon visqueuses}, 
  J. Math. Pures Appl. \textbf{76} (1997), 609--647. 

\bibitem{Da2}
\newblock R. Danchin, 
\newblock  \emph{Persistance de structures g\'eom\'etriques et 
  limite non visqueuse pour les fluides incompressibles en 
  dimension quelconque}, 
  Bull. Soc. Math. France \textbf{127} (1999), 179--227.

\bibitem{De}
\newblock J.-M. Delort,
\newblock \emph{Existence de nappes de tourbillon en dimension deux},
J. Amer. Math. Soc. \textbf{4} (1991), , 553–-586.

\bibitem{De1} 
\newblock W. Deng, 
\newblock \emph{Resolvent estimates for a two-dimensional 
non-selfadjoint operator}, preprint, 2010. 

\bibitem{De2} 
\newblock W. Deng, 
\newblock \emph{Pseudospectrum for Oseen vortices operators}, 
preprint, 2011. 

\bibitem{DR}
\newblock P. G. Drazin and W. H. Reid,
\newblock ``Hydrodynamic stability'', Second edition, Cambridge 
University Press, Cambridge, 2004.

\bibitem{FK}
\newblock H. Fujita and T. Kato,
\newblock \emph{On the Navier-Stokes initial value problem. I.},
Arch. Rational Mech. Anal. \textbf{16} (1964), 269–-315.

\bibitem{GG}
\newblock I. Gallagher and Th. Gallay,
\newblock \emph{Uniqueness for the two-dimensional Navier-Stokes 
  equation with a measure as initial vorticity}, 
Math. Ann., \textbf{332} (2005), 287--327.

\bibitem{GGN}
\newblock I. Gallagher, Th. Gallay and F. Nier,
\newblock \emph{ Spectral asymptotics for large skew-symmetric 
perturbations of the harmonic oscillator}, 
Int. Math. Res. Notices, \textbf{2009} (2009), 2147--2199.

\bibitem{GGL}
\newblock I. Gallagher, Th. Gallay and P.-L. Lions,
\newblock \emph{On the uniqueness of the solution of the two-dimensional
  Navier-Stokes equation with a Dirac mass as initial vorticity},
Math. Nachr., \textbf{278} (2005), 1665--1672.

\bibitem{Ga}
\newblock Th. Gallay,
\newblock \emph{Interaction of vortices in weakly viscous planar flows}, 
Arch. Ration. Mech. Anal., \textbf{200} (2011), 445–-490.

\bibitem{GR}
\newblock Th. Gallay and L. M. Rodrigues, 
\newblock \emph{Sur le temps de vie de la turbulence bidimensionnelle}, 
Ann. Fac. Sci. Toulouse Math. \textbf{17} (2008), 719–-733.

\bibitem{GW1}
\newblock Th. Gallay and C. E. Wayne,
\newblock \emph{Invariant manifolds and the long-time asymptotics 
of the Navier-Stokes and vorticity equations on {$\R^2$}},
Arch. Ration. Mech. Anal., \textbf{163} (2002), 209--258. 

\bibitem{GW2}
\newblock Th. Gallay and C.E. Wayne,
\newblock  \emph{Global stability of vortex solutions of the two-dimensional 
  Navier-Stokes equation}, 
Commun. Math. Phys., \textbf{255} (2005), 97--129.

\bibitem{GW3}
\newblock Th. Gallay and C.E. Wayne,
\newblock  \emph{Existence and stability of asymmetric Burgers 
vortices}, J. Math. Fluid Mech. \textbf{9} (2007), 243--261.

\bibitem{GK}
\newblock Y. Giga and T. Kambe,  
\newblock \emph{Large time behavior of the vorticity of two dimensional 
viscous flow and its application to vortex formation},  
Commun. Math. Phys., \textbf{117}, (1988) 549--568.

\bibitem{GMO}
\newblock Y.~Giga, T.~Miyakawa and H.~Osada, 
\newblock \emph{Two-dimensional Navier-Stokes flow with measures 
 as initial vorticity},
Arch. Rational Mech. Anal., \textbf{104} (1988), 223--250. 

\bibitem{Gr}
\newblock E. Grenier,
\newblock \emph{On the nonlinear instability of Euler and 
  Prandtl equations},
 Comm. Pure Appl. Math. \textbf{53} (2000), 1067--1091.

\bibitem{He} 
\newblock H. von Helmholtz,
\newblock \emph{\"Uber Integrale des hydrodynamischen Gleichungen, 
welche die Wirbelbewegungen entsprechen}; 
J. reine angew. Math., \textbf{55} (1858), 25-55.  

\bibitem{Hm1}
\newblock T. Hmidi,
\newblock \emph{R\'egularit\'e h\"old\'erienne des poches de 
  tourbillon visqueuses}, 
  J. Math. Pures Appl. \textbf{84} (2005), 1455--1495.

\bibitem{Hm2}
\newblock T. Hmidi,
\newblock \emph{Poches de tourbillon singuli\`eres dans un fluide 
  faiblement visqueux}, 
  Rev. Mat. Iberoamericana \textbf{22} (2006), 489--543. 

\bibitem{JMV}
\newblock J. Jim{\'e}nez, H. K. Moffatt and C. Vasco,
\newblock \emph{The structure of the vortices in freely decaying 
two-dimensional turbulence}, 
J. Fluid Mech. \textbf{313} (1996), 209--222.

\bibitem{JKN}
\newblock F. Jing, E. Kanso, and  P. Newton, 
\newblock \emph{Viscous evolution of point vortex equilibria: The 
collinear state}, Phys. Fluids \textbf{22} (2010), 123102.

\bibitem{Ka1}
\newblock T. Kato,
\newblock ``Perturbation theory for linear operators'', 
Grundlehren der mathematischen Wissenschaften \textbf{132}, 
Springer, New York, 1966. 

\bibitem{Ka3}
\newblock T. Kato,
\newblock \emph{Nonstationary flows of viscous and ideal fluids 
  in $\R^3$}, 
  J. Functional Analysis \textbf{9} (1972), 296--305.

\bibitem{Ka2}
\newblock T.~Kato,
\newblock \emph{The Navier-Stokes equation for an incompressible fluid 
  in $\R^2$ with a measure as the initial vorticity},
Differential Integral Equations, \textbf{7} (1994), 949--966. 

\bibitem{Ki}
\newblock G. R. Kirchhoff,
\newblock ``Vorlesungen \"uber Mathematische Physik. Mekanik.'',  
Teubner, Leipzig, 1876. 

\bibitem{DV} 
\newblock S. Le Diz\`es and A. Verga,
\newblock \emph{Viscous interactions of two co-rotating vortices 
before merging}, 
J. Fluid Mech., \textbf{467} (2002), 389--410. 

\bibitem{Le1}
\newblock J. Leray, 
\newblock \emph{\'Etude de diverses \'equations int\'egrales non 
lin\'eaires et de quelques probl\`emes que pose l'hydrodynamique},
J. Math. Pures Appl., IX. S\'er. 12 (1933), 1-82. 

\bibitem{Le2}
\newblock J. Leray, 
\newblock \emph{Sur le mouvement d'un liquide visqueux emplissant 
l'espace}, Acta Math. \textbf{63} (1934), 193–-248.

\bibitem{LL} 
\newblock E. Lieb and M. Loss, 
\newblock ``Analysis'', Graduate Studies in Mathematics, \textbf{14}, 
AMS, Providence, 1997.

\bibitem{McW}
\newblock J.C. Mc Williams,
\newblock \emph{The vortices of two-dimensional turbulence.}, 
J. Fluid. Mech., \textbf{219} (1990), 361--385. 

\bibitem{Ma} 
\newblock Y.~Maekawa,
\newblock \newblock{Spectral properties of the linearization at 
the Burgers vortex in the high rotation limit}, 
J. Math. Fluid Mech., \textbf{13} (2011), 515--532. 

\bibitem{Maj}
\newblock A. Majda
\newblock \emph{Remarks on weak solutions for vortex sheets with 
  a distinguished sign},
 Indiana Univ. Math. J. \textbf{42} (1993), 921--939.

\bibitem{MB} 
\newblock A. Majda and A. Bertozzi, 
\newblock "Vorticity and Incompressible Flow", 
Cambridge Texts in Applied Mathematics \textbf{27},  
Cambridge University Press, 2002.

\bibitem{Mar0}
\newblock C. Marchioro,
\newblock \emph{Euler evolution for singular initial data and vortex 
  theory: a global solution}, 
  Commun. Math. Phys. \textbf{116}(1988), 45--55.

\bibitem{Mar1}
\newblock C. Marchioro,
\newblock \emph{On the vanishing viscosity limit for two-dimensional 
Navier-Stokes equations with singular initial data}, 
Math. Methods Appl. Sci., \textbf{12} (1990), 463--470. 

\bibitem{Mar2}
\newblock C. Marchioro, 
\newblock \emph{On the inviscid limit for a fluid with a concentrated 
vorticity}, 
Commun. Math. Phys., \textbf{196} (1998), 53--65. 

\bibitem{MP1} 
\newblock C. Marchioro and M. Pulvirenti, 
\newblock \emph{Vortices and localization in Euler flows}, 
  Commun. Math. Phys. \textbf{154} (1993), 49--61. 

\bibitem{MP2} 
\newblock C. Marchioro and M. Pulvirenti, 
\newblock  ``Mathematical theory of incompressible nonviscous fluids'',  
  Applied Mathematical Sciences \textbf{96}, Springer, New York, 1994.

\bibitem{Mas} 
\newblock N. Masmoudi,
\newblock \emph{Remarks about the inviscid limit of the 
  Navier-Stokes system},  
  Comm. Math. Phys.  \textbf{270} (2007), 777--788.

\bibitem{MDL} 
\newblock P. Meunier, S. Le Diz\`es and T. Leweke,
\newblock \emph{Physics of vortex merging},  
Comptes Rendus Physique \textbf{6} (2005), 431--450. 

\bibitem{MKO} 
\newblock H. K. Moffatt, S. Kida and K. Ohkitani,
\newblock \emph{Stretched vortices---the sinews of turbulence; 
large-Reynolds-number asymptotics}, 
J. Fluid Mech., \textbf{259} (1994), 241--264. 

\bibitem{NSUW}
\newblock R. Nagem, G. Sandri, D. Uminsky and C.E. Wayne, 
\newblock \emph{Generalized Helmholtz-Kirchhoff model for two-dimensional 
 distributed vortex motion}, 
SIAM J. Appl. Dyn. Syst., \textbf{8} (2009), 160--179. 

\bibitem{Ne}
\newblock P. Newton, 
\newblock ``The $N$-vortex problem. Analytical techniques'', 
  Applied Mathematical Sciences \textbf{145}, Springer, New York, 
  2001. 

\bibitem{Os}
\newblock H. Osada, 
\newblock \emph{Diffusion processes with generators of generalized 
divergence form}, 
J. Math. Kyoto Univ. \textbf{27} (1987) 597–-619.

\bibitem{Pa}
\newblock A. Pazy,
\newblock ``Semigroups of linear operators and applications to partial 
differential equations'', Applied Mathematical Sciences \textbf{44}, 
Springer, New York, 1983.

\bibitem{PP}
\newblock A. Prochazka and D. I. Pullin,
\newblock \emph{On the two-dimensional stability of the axisymmetric 
  Burgers vortex},
Phys. Fluids \textbf{7} (1995), 1788--1790. 

\bibitem{RS}
\newblock M. Reed and B. Simon, 
\newblock ``Methods of modern mathematical physics. IV. Analysis of 
operators'',  Academic Press, New York, 1978.

\bibitem{Ru}
\newblock W. Rudin,
\newblock ``Real and complex analysis'', McGraw-Hill, New York, 1966. 

\bibitem{RWG}
\newblock M. A. Rutgers, X-l. Wu, and W. I. Goldburg, 
\newblock \emph{The Onset of 2-D Grid Generated Turbulence in 
Flowing Soap Films}, 
Phys. Fluids \textbf{8} (1996). 

\bibitem{SC1}
\newblock M. Sammartino and R. Caflisch,
\newblock \emph{Zero viscosity limit for analytic solutions 
  of the Navier-Stokes equation on a half-space. I. Existence for 
  Euler and Prandtl equations},
  Comm. Math. Phys. \textbf{\bf 192} (1998), 433--461. 

\bibitem{SC2}
\newblock M. Sammartino and R. Caflisch,
\newblock \emph{Zero viscosity limit for analytic solutions 
  of the Navier-Stokes equation on a half-space. II. Construction 
  of the Navier-Stokes solution},
  Comm. Math. Phys. \textbf{192} (1998), 463--491. 

\bibitem{SY} 
\newblock G. Sell and Y. You, 
\newblock ``Dynamics of evolutionary equations'',  
Applied Mathematical Sciences \textbf{143}, Springer, New York, 2002.

\bibitem{St} 
\newblock E. Stein,
\newblock ``Harmonic analysis: real-variable methods, orthogonality, 
and oscillatory integrals'', Princeton Mathematical Series, \textbf{43}, 
Princeton, 1993.

\bibitem{Su}
\newblock F. Sueur,
\newblock \emph{Vorticity internal transition layers for the 
  Navier-Stokes equations}, preprint, 2008. 

\bibitem{Sw} 
\newblock H. Swann, 
\newblock \emph{The convergence with vanishing viscosity of 
  nonstationary Navier-Stokes flow to ideal flow in $\R^3$}, 
  Trans. Amer. Math. Soc. \textbf{157} (1971), 373--397. 

\bibitem{Tao}
\newblock T. Tao, 
\newblock \emph{Localisation and compactness properties of the 
Navier-Stokes global regularity problem}, preprint, 2011. 

\bibitem{TK} 
\newblock L. Ting and R. Klein, 
\newblock ``Viscous vortical flows'', 
Lecture Notes in Physics \textbf{374}, Springer-Verlag, Berlin, 1991.

\bibitem{TT} 
\newblock L. Ting and C. Tung,  
\newblock \emph{Motion and decay of a vortex in a nonuniform stream},
Phys. Fluids, \textbf{8} (1965), 1039--1051.

\bibitem{TE} 
\newblock L. N. Trefethen and M. Embree, 
\newblock ``Spectra and pseudospectra: the behavior of nonnormal 
matrices and operators", Princeton University Press, 2005. 

\bibitem{Vil}
\newblock C.~Villani,
\newblock \emph{A review of mathematical topics in collisional kinetic 
  theory}, Handbook of mathematical fluid dynamics Vol. I, 
  pages 71--305, North-Holland, Amsterdam, 2002.

\bibitem{Vi} 
\newblock M. Vishik,
\newblock \emph{Incompressible flows of an ideal fluid with vorticity 
  in borderline spaces of Besov type}, 
  Ann. Sci. \'Ecole Norm. Sup. \textbf{32} (1999), 769--812.

\bibitem{Yu1}
\newblock V. Yudovich,
\newblock \emph{Non-stationary flows of an ideal incompressible 
  fluid. (Russian)},
  \v Z. Vy\v cisl. Mat. i Mat. Fiz. \textbf{3} (1963), 1032--1066.

\bibitem{Yu2}
\newblock V. Yudovich,
\newblock \emph{Uniqueness theorem for the basic nonstationary problem 
  in the dynamics of an ideal incompressible fluid},
  Math. Res. Lett. \textbf{2} (1995), 27--38. 

\end{thebibliography}
\end{document}